\documentclass{article}
\usepackage[english]{babel}
\usepackage[utf8]{inputenc}
\usepackage{amsmath,amsfonts,amsthm,amssymb}
\usepackage[alphabetic]{amsrefs}
\usepackage[all,cmtip]{xy}
\usepackage{graphicx}
\usepackage{verbatim}

\long\def\symbolfootnote[#1]#2{\begingroup%
\def\thefootnote{\fnsymbol{footnote}}\footnote[#1]{#2}\endgroup}

\newtheorem{theorem}{Theorem}[section]
\newtheorem{proposition}[theorem]{Proposition}
\newtheorem{corollary}[theorem]{Corollary}
\newtheorem{lemma}[theorem]{Lemma}

\newtheorem{question}[theorem]{Question}

\theoremstyle{definition}

\newtheorem{remark}[theorem]{Remark}
\newtheorem{definition}[theorem]{Definition}

\renewcommand{\proof}{\medskip\par\noindent\textbf{Proof:} \ignorespaces}
\renewcommand{\qed}{\quad\hskip0pt\null\hfill$\square$\par}

\newcommand*{\longhookrightarrow}{\ensuremath{\lhook\joinrel\relbar\joinrel\rightarrow}}
\newcommand{\df}[1]{{\bf{#1}}}
\newcommand{\Ends}{{\text{Ends}}}

\def\F{\mathcal{F}}
\def\E{\mathbb{E}}

\def\Z{\mathbb{Z}}
\def\R{\mathbb{R}}
\def\N{\mathbb{N}}
\def\H{\mathbb{H}}
\def\O{\mathcal{O}}

\def\cD{\mathcal{D}}
\def\C{\mathbb{C}}

\newcommand{\Cayley}{\mathrm{Cayley}}

\newcommand{\Isom}{\text{Isom}}

\newcommand{\mon}{{\text{Mon}}}
\newcommand{\aut}{{\text{Aut}}}
\newcommand{\Stab}{{\text{Stab}}}

\title{The topology of the minimal regular cover of the Archimedean tessellations}

\begin{document}
\maketitle

\begin{center}\bf
Thierry Coulbois$^a$, Daniel Pellicer$^b$\symbolfootnote[2]{Partially supported by PAPIIT-Mexico under grant IB100312-2.}, Miguel Raggi $^c$\symbolfootnote[3]{Partially supported by DGAPA Posdoctoral scholarship and PAPIIT-Mexico under grant IB100312-2.}, Camilo Ram\'irez$^d$\symbolfootnote[4]{Partially supported by CONACYT and CCM-UNAM, Morelia.}  \& Ferr\'an Valdez$^e$\symbolfootnote[5]{Partially supported by LAISLA, CONACYT CB-2009-01 127991 and PAPIIT projects IN103411 \& IB100212.}
\end{center}

\begin{center}\it
$^a$ LATP, Universit\'e d'Aix-Marseille \\
39, rue F. Joliot Curie\\
13453 Marseille, Cedex 13, France\\
\emph{e-mail: }\texttt{thierry.coulbois@univ-amu.fr}
\end{center}
\begin{center}\it
$^{b,}$ $^{c,}$ $^{d,}$ $^{e}$ 
Centro de Ciencias Matem\'aticas,
UNAM, Campus Morelia \\
C.P. 58190, Morelia, Michoac\'an,  M\'exico.\\
$^{b}$\emph{e-mail: }\texttt{pellicer@matmor.unam.mx}\\
$^{c}$\emph{e-mail: }\texttt{mraggi@gmail.com}\\
$^{d}$\emph{e-mail: }\texttt{camilomaluendas@matmor.unam.mx}\\
$^{e}$\emph{e-mail: }\texttt{ferran@matmor.unam.mx}\\
\end{center}

\begin{abstract}
\noindent In this article we determine, for an infinite family of maps on the plane, the topology of the surface on which the minimal regular covering occurs. This infinite family includes all Archimedean maps.
\end{abstract}

\section{Introduction}

A map on a surface $S$ follows the idea of the map of the Earth. The surface is decomposed into ``countries'' (faces) where every border (edge) belongs to precisely two countries. The points where three or more countries meet correspond to the vertices of the map.

Symmetries of maps on surfaces have attracted increasing attention in the last decades and most of the attention has been devoted to regular maps on compact surfaces \cite{conderdob}, \cite{Gareth}, \cite{Wilson}.

In order to work with less symmetric maps, in \cite{Hartley} the author provides a procedure to describe an arbitrary map $M$ algebraically in terms of a regular cover $\widetilde{M}$. In particular, Hartley showed a way of representing $M$ as a pair consisting of $\widetilde{M}$ and a subgroup of its automorphisms. Hartley did this work in the context of abstract regular polytopes, however, the ideas and techniques apply directly for maps due to the similarity of their structure with that of rank 3 abstract polytopes.

A regular cover $\widetilde{M}$ of the map $M$ can be infinite (have infinitely many vertices, edges and faces) even when $M$ is finite, and as a consequence Hartley's technique is hard to apply when attacking some concrete problems. In those instances it may be convenient to require $\widetilde{M}$ to be as small as possible. In particular, if $M$ is finite, we can require $\widetilde{M}$ to be finite as well. Such a regular map always exists and is called the \df{minimal regular cover}.

Hartley's procedure can be applied to maps on either compact or non-compact surfaces. In order to better understand his technique, subsequent papers determined the minimal regular cover for some relevant families of maps, like those arising from the Archimedean solids \cite{HarWil} and those arising from the prisms and antiprisms \cite{PelGor}. Since all these maps have a finite number of vertices, edges and faces, the minimal regular covers also have a finite number of vertices, edges and faces, and lie on compact surfaces which can be easily derived from the orientability and the Euler characteristic.

The Archimedean tessellations of the plane were the first maps with infinitely many vertices for which the minimal regular cover was determined (see \cite{GorDan}, \cite{Mixer}), however, determining the underlying surface of these covers proved more challenging and remained open.

In this paper we show that the minimal regular cover of any map on a large family of maps on the Euclidean plane, including all Archimedean tessellations, lies on a topological surface obtained by glueing infinitely many torii along a ray (see figure \ref{F:1}). This surface is also known as the \df{Loch Ness Monster}. This nomenclature is due to Ghys \cite{Ghys}.
\begin{center}
    \label{F:1}
\includegraphics[scale=0.8]{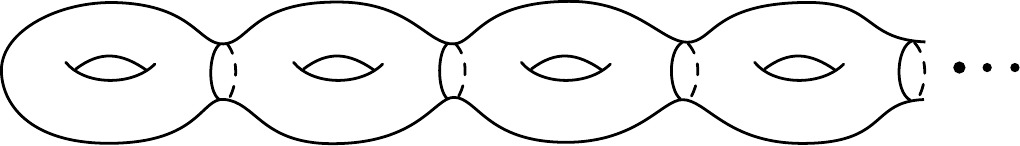}\\
\small{Figure 1: The Loch Ness monster.}\\
\end{center}
The \textit{Loch Ness monster} is an ubiquitous surface: it appears as the generic leaf in orientable laminations of compact spaces \cite{Ghys}, also as the generic leaf of a homogeneous holomorphic foliation in $\C^2$ \cite{V1}, as the flat surface associated to the billiard game on a generic polygon \cite{V2}, and even when trying to realize any countable subgroup of ${\rm GL}_+(2,\R)$ as the Veech group of a flat surface \cite{PSV}. As we will see, its appearance in the context of Euclidean tessellations is a consequence of a classical result due to Stallings.

\section{Tessellations and covers}
  \label{S:tess}

Throughout we consider a \df{map} to be a 2-cell embedding of a \df{locally finite} graph (graph whose vertices belong to finitely many edges) on a surface. In other words, the vertices of a graph are mapped bijectively to a discrete set of points of the surface and the edges are mapped to arcs between the corresponding points in such a way that the the following properties hold.
\begin{itemize}
  \item No edge self-intersects.
  \item The images of two edges intersect if and only if the edges have a vertex in common. In this case the images intersect only in the image of the vertex in common.
  \item If we remove the image of the graph, the remaining connected components (called \df{faces}) are homeomorphic to discs.
\end{itemize}
Hence, a map consists of its graph, its surface, and an embedding. We shall abuse notation and interpret the vertices and edges of the graph as their images on the surface. For simplicity, we shall denote the map $\Gamma\stackrel{i}{\hookrightarrow} S$, where the graph $\Gamma$ is embedded to the surface $S$ by the mapping $i$, simply by $M=M(\Gamma,i,S)$.

Two maps $M_1$ and $M_2$ on surfaces $S_1$ and $S_2$ are \df{isomorphic} whenever there is an homeomorphism $\phi$ from $S_1$ to $S_2$ that induces an isomorphism between the graphs of $M_1$ and $M_2$. Thus, the isomorphism between $S_1$ and $S_2$ consists of $\phi$ together with the isomorphism of the graphs of $M_1$ and $M_2$. We are interested on maps modulo isomorphisms. An \df{automorphism} of a map is an isomorphism of the map to itself. Since we consider isomorphic maps as equal, an automorphism of a map can be alternatively defined as an automorphism of its graph that can be extended to an homeomorphism of the surface to itself. We denote the group of automorphisms of a map $M$ by $\aut(M)$.

Whenever we consider a map $M$ on a surface $S$ with a specified metric we say that an isometry $\alpha$ of $S$ is a \df{symmetry} of $M$ whenever $\alpha$ preserves the vertex and edge sets of $M$ setwise. Note that every symmetry induces an automorphism of the map, but the converse is not always true. In this paper we shall devote particular attention to maps on the Euclidean plane $\mathbb{E}^2$, where the symmetries are with respect to the standard metric.

A \df{flag} of a map is a triangle on the surface whose vertices are a vertex $v$ of the map, the midpoint of an edge $e$ containing $v$, and an interior point of a face containing $e$. We may assume that for all flags containing the face $F$, the same interior point of $F$ is chosen to be a vertex of the corresponding triangle. By doing this, every map induces the triangulation given by its flags of the underlying surface. In this paper we are interested only on maps whose graph does not have loops, with the property that every edge belongs to two distinct faces; in this case every flag can be identified with an ordered incident triple of vertex, edge and face. Henceforth we shall abuse notation and understand flags either as triangles or as ordered triples. We shall denote the set of orbits of flags under $\aut(M)$ of the map $M$ by $\O_{\aut}(M)$. In figure \ref{F:2} we find an example of a cube divided into flags. There is a marked flag $\Phi$ with its corresponding vertex $v$, edge $e$ and face $f$. Given a flag $\Phi$, there is a unique flag $\Phi^0$ (resp. $\Phi^1$ and $\Phi^2$) that differs from $\Phi$ precisely on the vertex (resp. on the edge and on the face). Then, if $w$ is a word on the set $\{0, 1, 2\}$ and $j \in \{0,1,2\}$, we define recursively $\Phi^{wj} := (\Phi^w)^j$. The flag $\Phi^i$ is called the \df{$i$-adjacent flag of $\Phi$}. In figure \ref{F:2}, we find an example of the cube with some flags marked with their name.
\begin{center}
    \label{F:2}
  \includegraphics[scale=0.6]{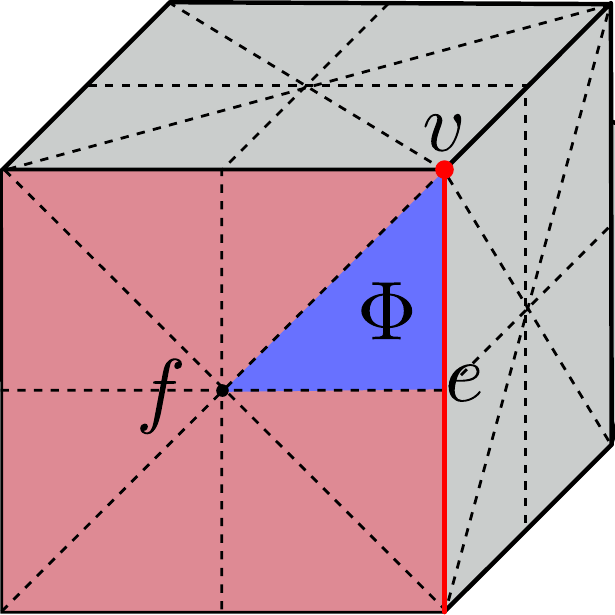}\hspace{6mm}\includegraphics[scale=0.6]{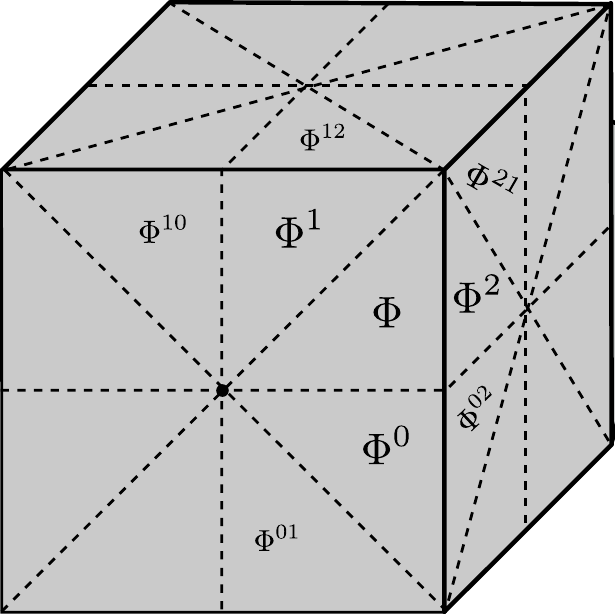}\\
    \small{Figure 2: A cube divided into flags.}
\end{center}

It follows directly from the definition that every automorphism $\varphi$ of a map $M$ is completely determined by the image of any of the flags of $M$ under $\varphi$. In other words, for every pair $\Phi$, $\Psi$ of flags of $M$ there is at most one automorphism of $M$ mapping $\Phi$ to $\Psi$. Hence, we regard as the most symmetric maps those in which the automorphism group acts transitively (regularly) on its flags. Maps with this property are called \df{regular}.

Note that if $M$ is a regular map on a surface $S$ then any flag can be chosen as a fundamental domain for $\aut(M)$ implying that {\em $\aut(M)$ acts on $S$ cocompactly}. Furthermore, given a point $p$ on $S$, the orbit of $p$ under $\aut(M)$ contains precisely one point on each flag of $M$, implying that $\aut(M)$ acts properly discontinuously on $S$. Recall that an action of a group $G$ on a topological space $X$ is \textbf{properly discontinuous} if for any compact set $K \subset X$, there are finitely many elements $g$ of $G$ such that $K\cap gK\neq\emptyset$.

The automorphism group of a regular map is generated by three involutions $\rho_0$, $\rho_1$ and $\rho_2$ which map a given base flag $\Phi$ to $\Phi^0$, $\Phi^1$ and $\Phi^2$ respectively. In fact, if $\varphi\in\aut(M)$ is such that $\varphi(\Phi) = \Phi^{j_1, \dots, j_k}$ for some $j_1, \dots, j_k \in \{0, 1, 2\}$, it can be proved with an inductive argument that $\varphi = \rho_{j_k} \rho_{j_{k-1}} \cdots \rho_{j_1}$ (see for example \cite[Propositions 2B7, 2B8]{ARP} with $J = \{0,1,2\}$).

The generators $\rho_i$ of the automorphism group of a regular map satisfy that $\rho_0 \rho_2 = \rho_2 \rho_0$ and that $\rho_0, \rho_1, \rho_2, \rho_0 \rho_1, \rho_0 \rho_2$, and $\rho_1 \rho_2$ are all different from $id$. Conversely, if $G$ is generated by involutions $\rho_0, \rho_1$ and $\rho_2$, with the properties that $\rho_0$ and $\rho_2$ commute and that the products of one or two distinct generators are different from $id$, then $G$ is the automorphism group of a unique regular map where the generator $\rho_i$ maps a base flag to its $i$-adjacent flag for $i\in \{0, 1,2\}$. This map can be constructed by considering a triangle with edges labelled $0, 1, 2$ for each element of $\langle \rho_0, \rho_1,\rho_2 \rangle$. The sides with label $i$ of the triangles $g_1$ and $g_2$ are identified whenever $g_1 = g_2 \rho_i$.

Examples of regular maps are the Platonic solids, viewed as maps on the sphere, the three regular tessellations $\{3,6\}$, $\{4,4\}$ and $\{6,3\}$ of the Euclidean plane and the regular tessellations $\{p,q\}$ with $1/p + 1/q < 1/2$ in the hyperbolic plane. Here $\{p,q\}$ indicates the tessellation obtained by arranging $q$ $p$-gonal faces around each vertex with no further identifications. It is well-known that the automorphism group of the tessellation $\{p,q\}$ is the Coxeter group
\[[p,q] = \langle \rho_0, \rho_1, \rho_2 \mid \rho_0^2 = \rho_1^2 = \rho_2^2 = (\rho_0 \rho_2)^2 = (\rho_0 \rho_1)^p = (\rho_1 \rho_2)^q = id\rangle.\]

Among the most studied infinite maps are the 11 tessellations of the Euclidean plane with (geometrically) regular convex polygons as faces and whose symmetry group acts transitively on the vertices. Three of them are the regular tessellations by equilateral triangles, by squares and by regular hexagons. The remaining eight are known as the \df{Archimedean tilings} of the plane.

Each Archimedean tiling is completely determined by the sequence of numbers of edges of the faces around a given vertex, given in the order in which they occur. We shall refer to each individual Archimedean tiling by the notation introduced by Gr\"unbaum. In this notation, the numbers of edges of faces are listed in the appropriate order and the Archimedean tilings are the ones listed below (see \cite[Section 2.1]{GruShe}).\\

\begin{tabular}{ccc}
\includegraphics[scale=0.24]{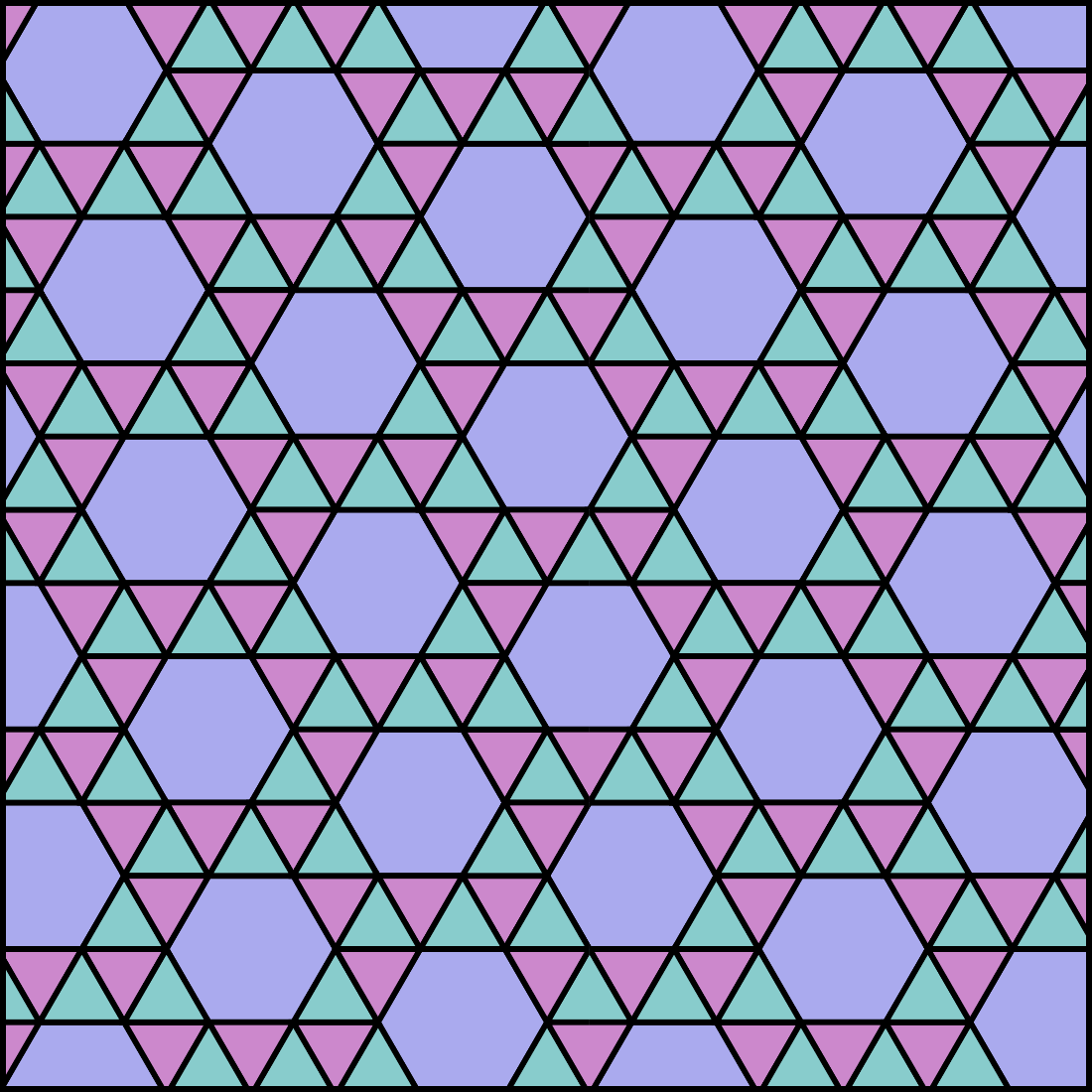} & \includegraphics[scale=0.24]{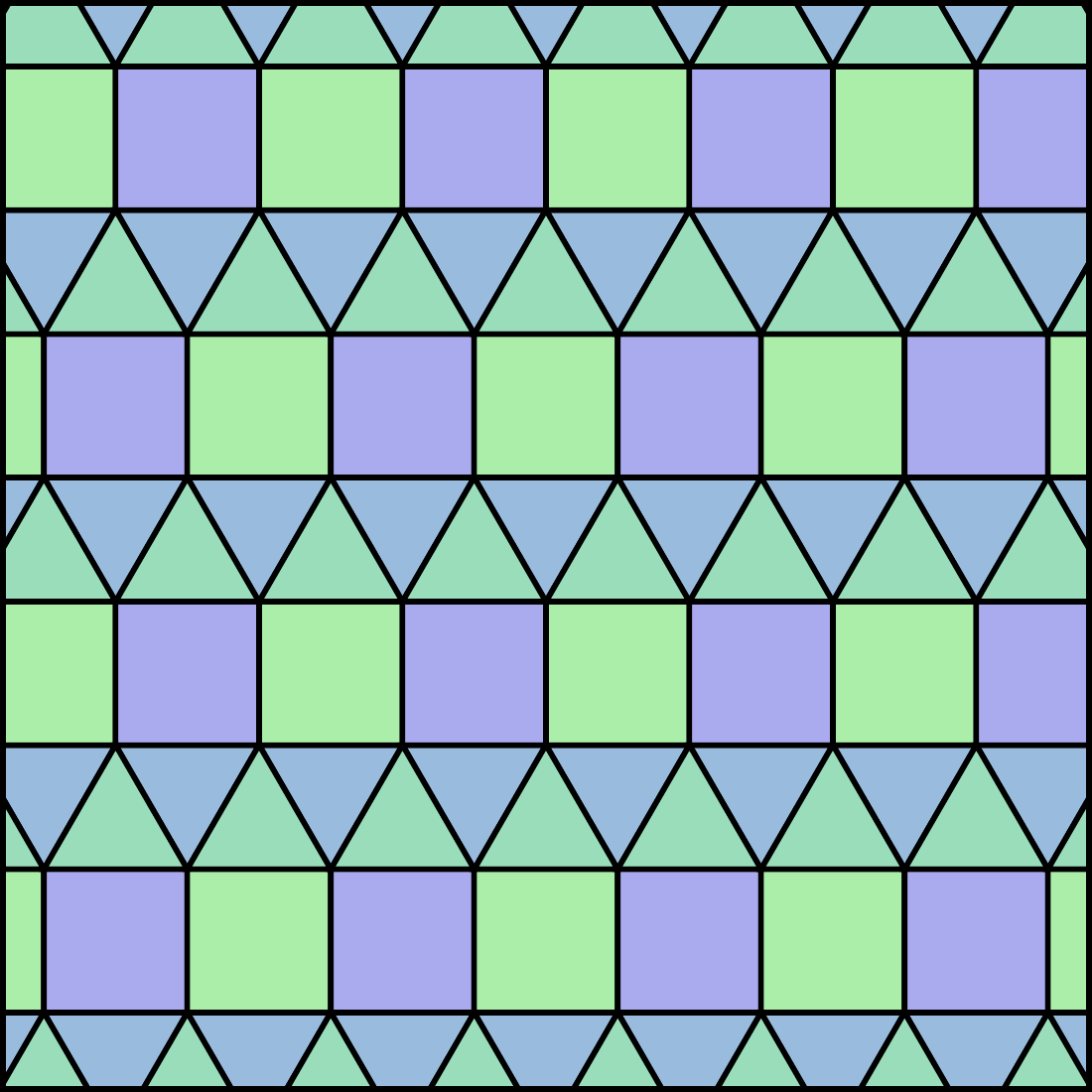} & \includegraphics[scale=0.24]{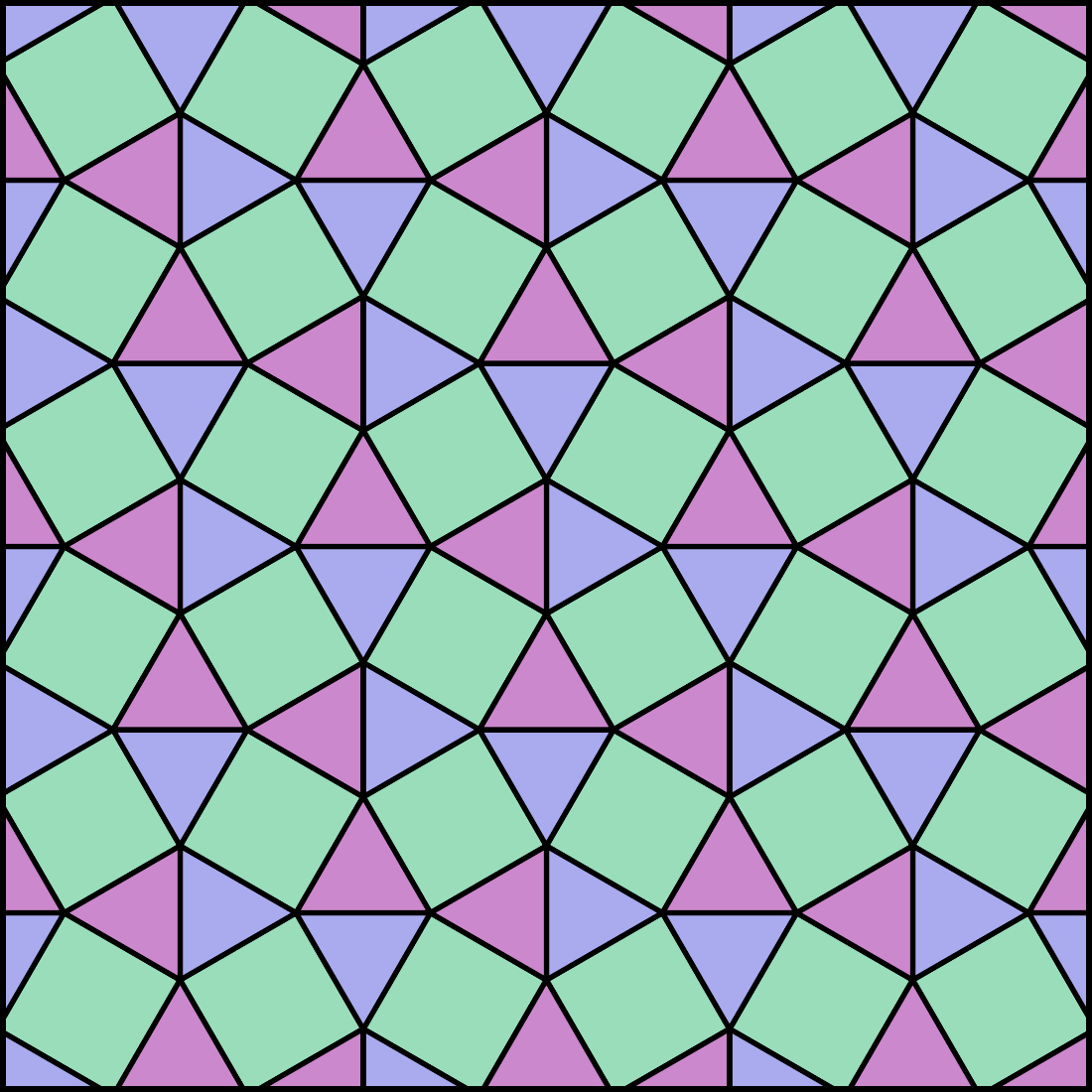} \\
3.3.3.3.6 & 3.3.3.4.4 & 3.3.4.3.4 \\

& & \\

\includegraphics[scale=0.24]{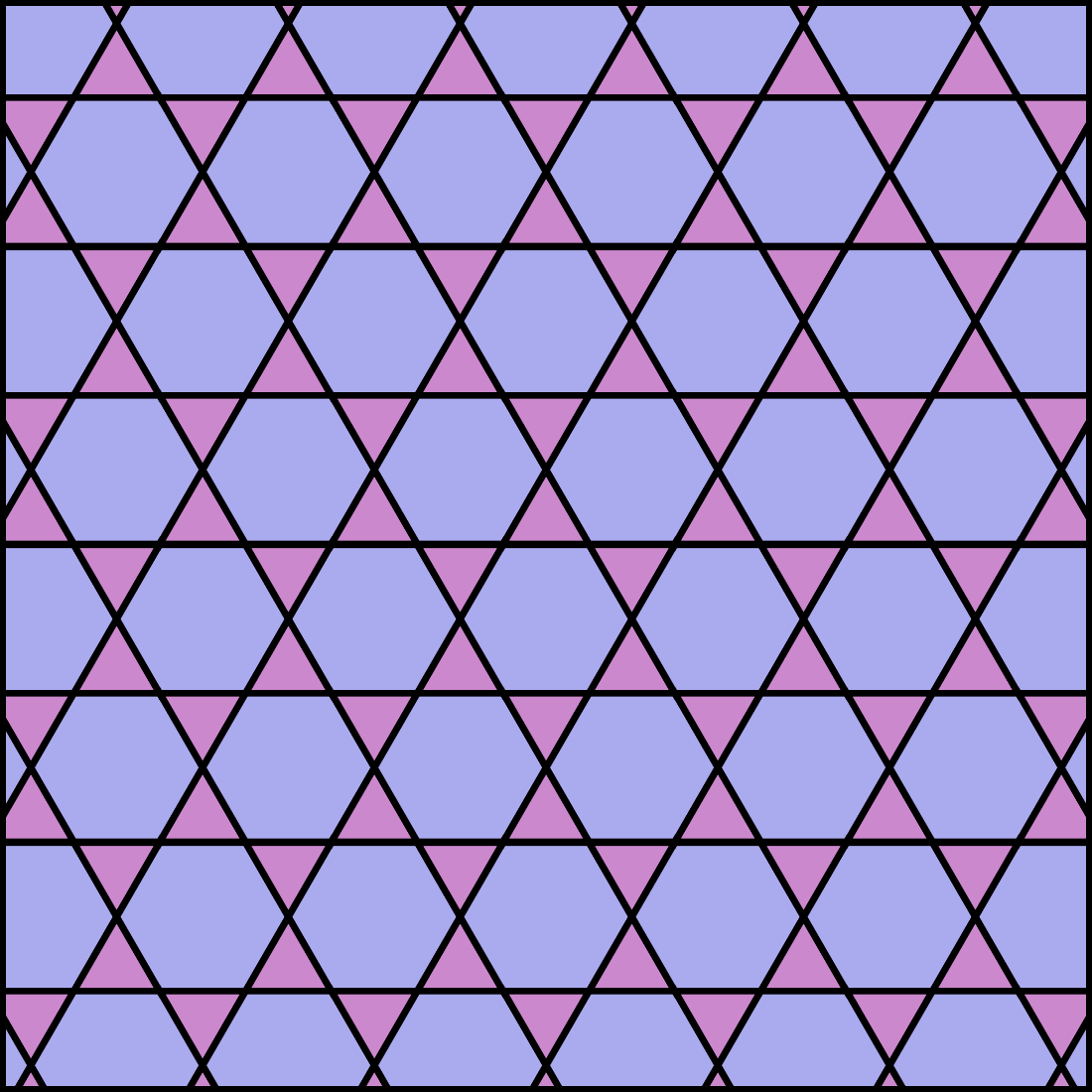} & \includegraphics[scale=0.24]{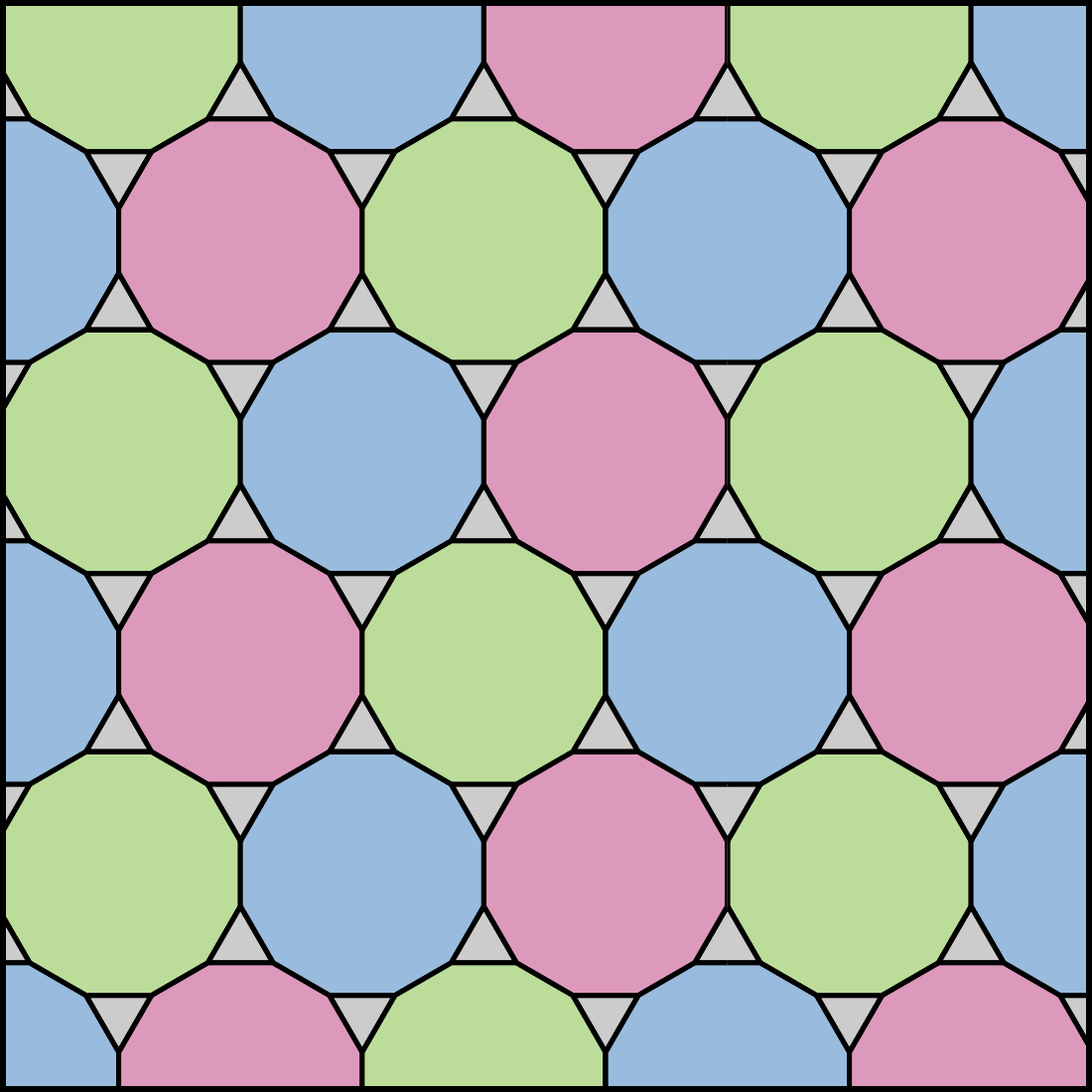} & \includegraphics[scale=0.24]{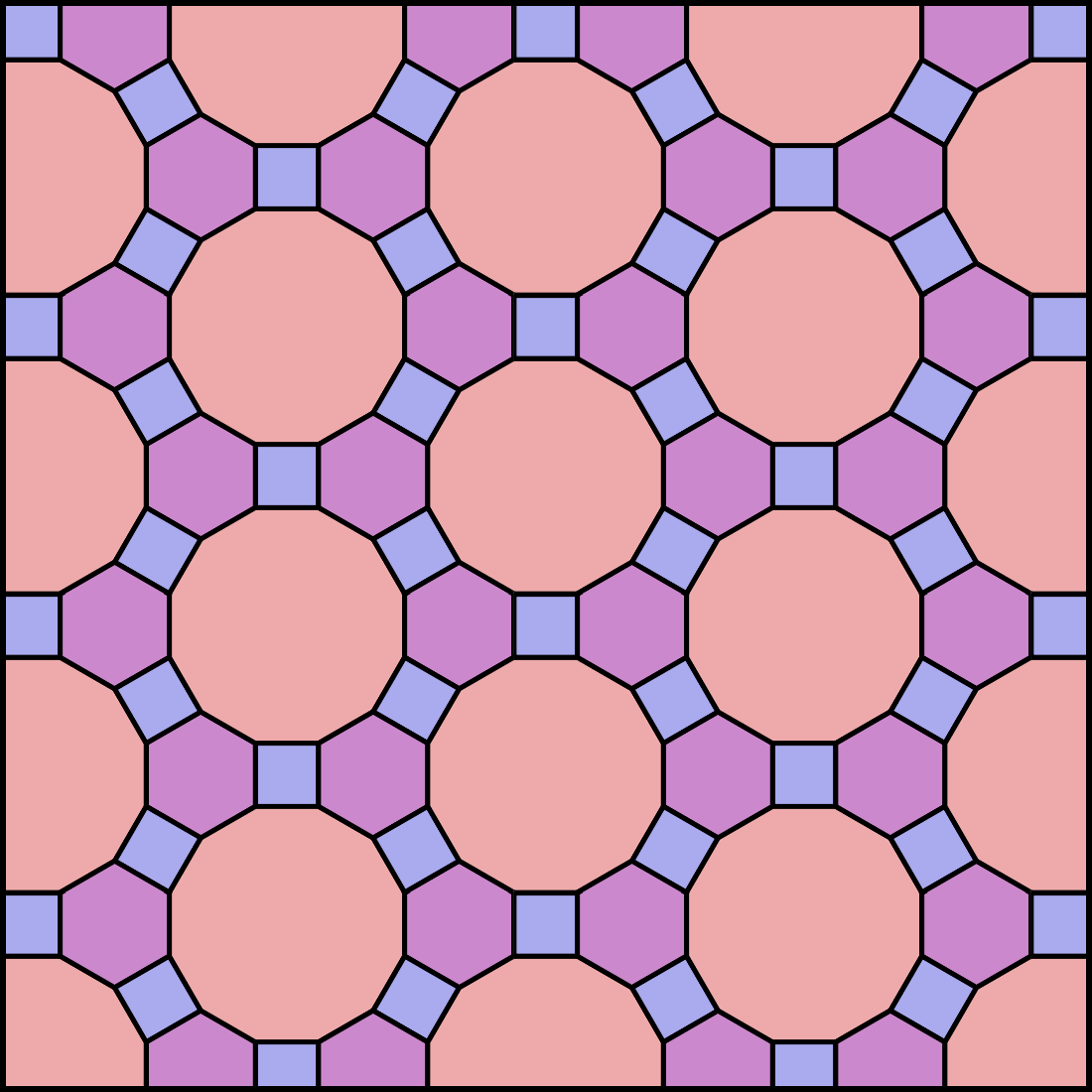} \\
3.6.3.6 & 3.12.12 & 4.6.12 \\

& & \\

\includegraphics[scale=0.24]{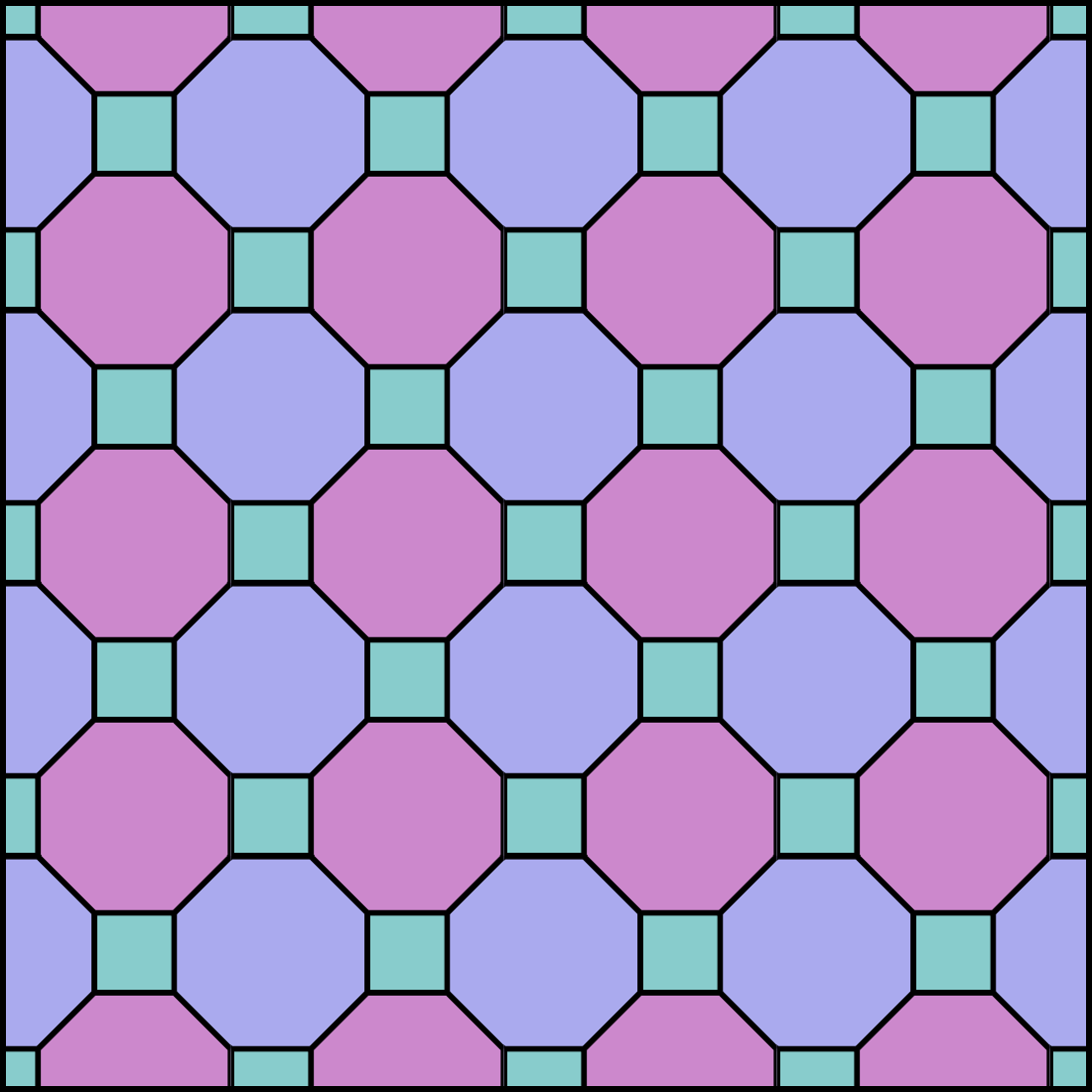} & \includegraphics[scale=0.24]{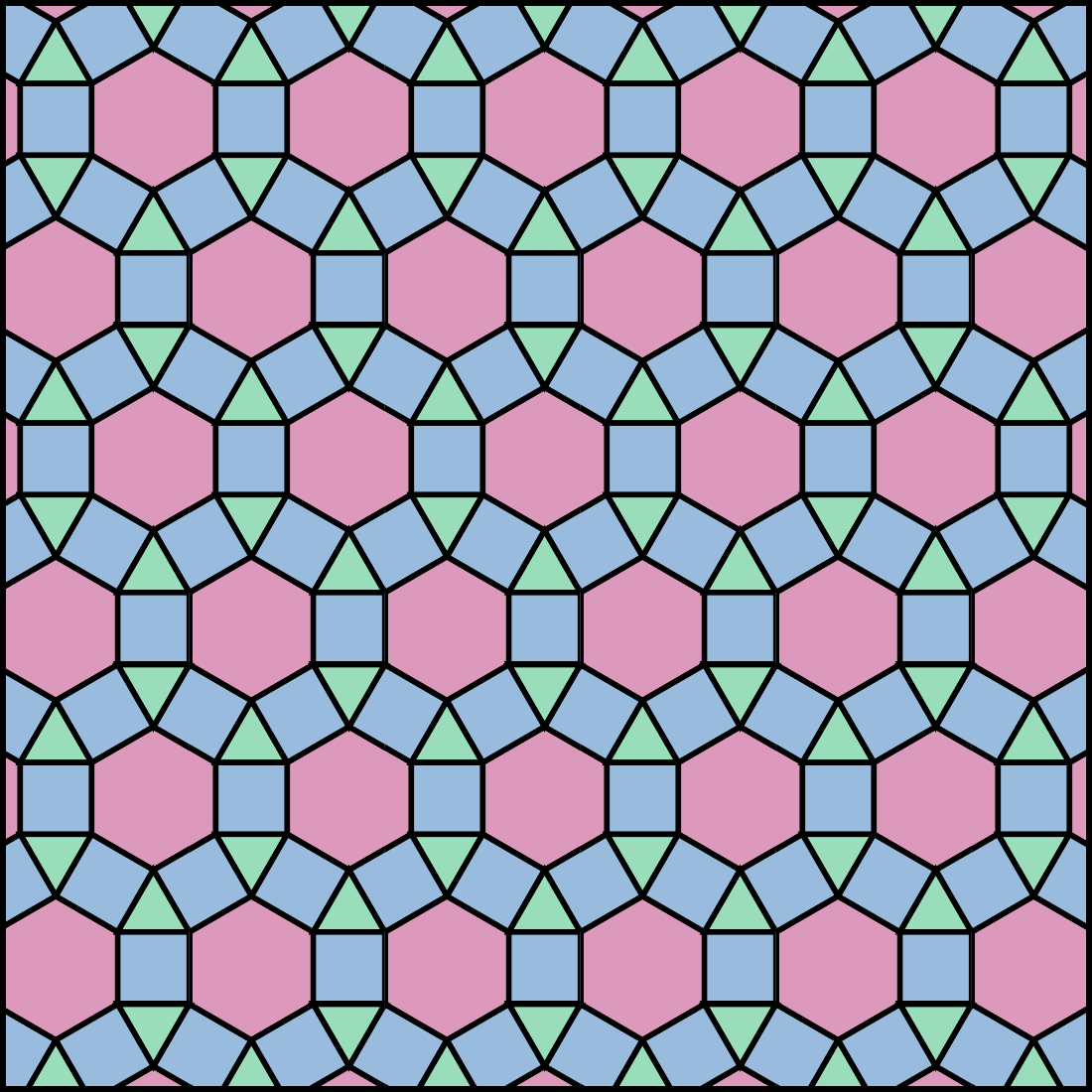} & \quad \\
4.8.8 & 3.4.6.4 & \\

\end{tabular}

\begin{center}
    \label{F:3}
    \small{Figure 3: The eight Archimedean tessellations. \\ Images by R. A. Nonenmacher, distributed under CC BY-SA 2.5.}
\end{center}


If $M$ is a map on the Euclidean plane whose group of isometries contains two translations with respect to linearly independent vectors, then $\aut(M)$ is one of the 17 crystallographic groups of $\mathbb{E}^2$. In particular, the translation subgroup of $\aut(M)$ has finite index in $\aut(M)$ (see \cite[Section 1.4]{GruShe}, \cite[Chapters 3--6]{Con} or \cite[Section 7.4]{Ratcliffe}). Furthermore, every parallelogram defined by the two generating translations intersects finitely many flags. The previous discussion proves the following proposition.
\begin{proposition}
\label{prop:finiteorbits}
Let $M$ be a map on the Euclidean plane whose group of isometries contains two translations with respect to linearly independent vectors. Then the translation subgroup has finite index in $\aut(M)$, and $\O_{\aut(M)}$ is finite.
\end{proposition}

Since, the symmetry group of a map $M$ can be embedded into $\aut(M)$, we have the following remark.

\begin{remark}
Let $M$ be a map as in Proposition \ref{prop:finiteorbits}, then the translation subgroup has finite index in the symmetry group of $M$.
\end{remark}

For $i \in \{0,1,2\}$ let $r_i$ be the involution on the set of flags $\F(M)$ of $M$ that maps every flag $\Phi$ to its $i$-adjacent. Then the \df{monodromy group} $\mon(M)$ of $M$ is the subgroup of the symmetric group on $\F(M)$ given by
$$\mon(M) = \langle r_0, r_1, r_2 \rangle.$$
Note that the involution $r_i$ needs not be a morphism of the map and needs not induce an homeomorphism of the surface; it is merely a bijection in the set of flags.
Clearly {\em the monodromy group of any map acts transitively on its flags}.

It follows from the previous definitions that for any flag $\Phi$ of a map $M$, any automorphism $\varphi$ of $M$ and $i \in \{0, 1, 2\}$,
$$\varphi[(\Phi)r_i] = [\varphi(\Phi)]r_i.$$
An inductive procedure implies that for any $\varphi \in\aut(M)$ and any $w \in\mon(M)$, $\varphi[(\Phi)w] = [\varphi(\Phi)]w$. Henceforth we shall assume a left action of $\aut(M)$ and a right action of $\mon(M)$ on $\F(M)$.

Whenever a map $M$ is regular there is a group isomorphism between $\aut(M)$ and $\mon(M)$ mapping $\rho_i$ to $r_i$ for $i \in \{0, 1, 2\}$. This can be proved inductively, but note that this isomorphism reverses the order on which the generators are applied to a given flag (see \cite[Theorem 3.9]{Monson}).

Given a map $M$, the generators $r_0, r_1, r_2$ of $\mon(M)$ satisfy that $r_0 r_2 = r_2 r_0$ and that the product of one or two distinct generators is different from $id$. Hence, there is a regular map $\widetilde{M}$ naturally associated to $M$ where $\aut(\widetilde{M})=\mon(M)$. If $M$ is regular then $M \cong \widetilde{M}$, otherwise $\widetilde{M}$ is a proper \df{regular cover} of $M$ (a covering map which is itself regular), possibly with branch points. (Here we think of coverings of surfaces mapping the vertex and edge set of one map onto the vertex and edge set of the other map.) Furthermore, $\widetilde{M}$ has the property of being covered by any regular cover of $M$. With this in mind, $\widetilde{M}$ is defined to be the \df{minimal regular cover} of $M$ (see \cite{Barry} and \cite{Monson} for further details).\\
\\
\indent We conclude this section with a technical lemma that will be used in \S \ref{S:main} to prove our main results.

%
\begin{lemma}
  \label{egamma}
Let $M(\Gamma,i,\E^{2})$ be a map whose faces are convex polygons satisfying the following.
\begin{itemize}
\item $\aut(M)$ contains two linearly independent translations $T_{1}, T_{2}\in\Isom(\E^{2})$.
\item The map $M$ has  at least two faces of different sizes.
\end{itemize}
Let $\pi : \widetilde{S} \rightarrow \E^{2}$ be the covering from the surface $\widetilde{S}$ of the minimal regular cover of $M$ to $\E^{2}$. Then there is a simple closed curve $\gamma$ in $\E^{2}$ with the following properties.
\begin{enumerate}
   \item The preimage of $\gamma$ under $\pi$ consists of disjoint union of simple closed curves.
   \item If $\widetilde{\gamma_1}$ is a simple closed curve in the preimage of $\gamma$, then the restriction of $\pi$ to $\widetilde{\gamma_1}$ is a bijection.
   \item There is a branch point of $\pi$ in the compact region bounded by $\gamma$.
\end{enumerate}
\end{lemma}

\proof Let $q$ be the least common multiple of the numbers of edges of all faces of $M$. Then all faces in the minimal regular cover of $M$ have precisely $q$ edges and the map $M$ has a branch point in the interior of each face with less than $q$ edges. Let $v_0$ be a branch point of $\pi$.\\ On the other hand, the map $M$ has two linearly independent translations $T_1, T_2$.
The curve $\gamma$ can be found as follows. Take a flag $\Phi$ in $M$. Pick a point $x\in\Phi$ \emph{which is not a branch point of $\pi$}. We can then find $k,l\in\N$ such that the parallelogram $P$ with vertices
$$
\{T_{1}^{k}(x),T_{2}^{l}T_{1}^{k}(x),T^{-k}_{1}T_{2}^{l}T_{1}^{k}(x),x=T^{-l}_{2}T^{-k}_{1}T_{2}^{l}T_{1}^{k}(x)\}
$$
contains $v_0$ in its interior. Clearly the choice of $x$ can be made such that $P$ contains no vertex of a flag of $M$ (that is, no vertex, midpoint of edge, or center of face of $M$). Let $w_1$ (resp. $w_2$) be the word on the generating set $\{r_0, r_1, r_2\}$ of $\mon(M)$ mapping $\Phi$ to $T_{1}^{k}(\Phi)$ (resp. $T_{2}^{l}(\Phi)$) induced by the crossings of one side of $P$ with the boundaries of the flags of $M$. Note that, starting on $\Phi$, no flag of $M$ is traversed more than once when applying $w_1$ or $w_2$. Up to raising both $w_{1}$ and $w_{2}$ to some power, we can suppose that each of these elements fixes every orbit in $\O_{\aut}(M)$ and acts on each flag of $\F(M)$ as a translation or as the identity. If $P$ contains no branch point of $\pi$ then take $P=\gamma$.
If $P$ contains branch points of $\pi$, we can perform a smooth homotopy on $P$ that leaves the vertices of the parallelogram fixed and avoids the branch points. (Note here that there can be no branch points on the edges of $M$ and there can be at most one branch point on the interior of each face of $M$. This implies that the set of branch points is discrete.)
\\
We claim that every lift of $P$ is a closed curve. In fact, the walk induced by $w^{-l}_{2}w^{-k}_{1}w_{2}^{l}w_{1}^{k}$ in the flags of $M$ lifts to a walk induced by the same word on $\mon(M)$ in the flags of its minimal regular cover. But $w^{-l}_{2}w^{-k}_{1}w_{2}^{l}w_{1}^{k}$ is trivial in $\mon(M)$ and $P$ contains no branch points itself implying that every lift of $P$ must be closed. Furthermore, if the restriction of $\pi$ to some lift of $P$ was a mapping $k$ to $1$, then there is a subword of $w^{-l}_{2}w^{-k}_{1}w_{2}^{l}w_{1}^{k}$ that fixes all flags of $M$. (Recall that a word $w$ in $\mon(M)$ induces a closed walk on the flags of the minimal regular cover of $M$ if and only if $w$ fixes all flags of $M$.) However this contradicts our choice of $w_1$ and $w_2$.\qed

\section{Ends and groups}
In this section we make a short review of the theory of ends of topological spaces and finitely generated groups.

\subsection{Ends of topological spaces} Although the discourse about ends will be presented in the most general context, we will apply it later to relatively nice topological spaces $X$: 2-dimensional real manifolds (\emph{i.e.} surfaces) and locally finite graphs.
Roughly speaking, an \emph{end} of a topological space is a ``point at infinity''. The set of ends forms a natural extension of $X$ and can be endowed with a topology. The (topological) space of ends is a topological invariant. As we will see later, the space of ends plus the genus, is \emph{the} topological invariant of an orientable surface. In the following paragraphs we present two equivalent definitions of the space of ends of a locally compact, locally connected, connected Hausdorff space $X$.
\begin{definition}\cite{Mil}
    \label{D:Milnor}
An \df{end} of $X$ is a function $E$ which assigns to each compact subset $K\subset X$ precisely one connected component $E(K)$ of the complement $X\setminus K$, subject to the requirement that $E(K)\supset E(L)$ whenever $K\subset L$. We denote the set of ends of $X$ by ${\rm Ends}(X)$.
\end{definition}
\begin{remark}
Remark that if $X$ is compact, then ${\rm Ends}(X)=\emptyset$.
\end{remark}
We define the topology in the space of ends as follows. Let $\mathcal{K}(X):=\{K_{i}\}_{i\in I}$ be the set of all compact subsets of $X$ indexed by some set $I$. We define the partial order $i\leq j$ if $K_{i}\subset K_{j}$. Define $A_{i}=\pi_{0}(X\setminus K_{i})$
, the set of connected components of $X\setminus K_{i}$. The set $A_{i}$ is always finite (see \cite[Lemma~1.1]{Ray}) and will be endowed henceforth with the discrete topology.
If $K_{i}\subset K_{j}$ then every connected component of $X\setminus K_j$ is contained in a unique connected component of $X\setminus K_i$. Hence, we have the natural maps
\begin{equation}
f_{ij}:A_{j}\to A_{i}
\end{equation}
that define an inverse system. Note that, by definition,
\begin{equation}
{\rm Ends}(X)=\varprojlim_{i\in I}A_{i}
\end{equation}
as sets. Henceforth, we will consider ${\rm Ends}(X)$ endowed with the limit topology. The space ${\rm Ends}(X)$ is compact, closed, has no interior points and is totally disconnected (see Theorem 1.5, \cite{Ray}).

We now present a second definition of $\mathrm{Ends}(X)$ due to Freudenthal \cite{F}. It is motivated by previous work of Caratheodory and Mazurkiewicz (see \cite{M}) and generalizes the notion of ``Randst\"ucke'' (literally ``boundary pieces'') introduced by B\'ela and Ker\'ekj\'art\'o for surfaces (see \cite{L} for a review) and later used by Richards \cite{R}.
\begin{definition}\cite{F}
Let $U_{1}\supseteq U_{2}\supseteq\ldots$ be an infinite sequence of non-empty connected open subsets of $X$ such that for each $i\in\N$ the boundary $\partial U_{i}$ is compact and $\bigcap\limits_{i\in\N}\overline{U_{i}}=\emptyset$. Two such sequences $U_{1}\supseteq U_{2}\supseteq\ldots$ and $U'_{1}\supseteq U'_{2}\supseteq\ldots$ are said to be equivalent if for every $i\in\N$ there exist $j$ such that $U_{i}\supseteq U'_{j}$ and $U'_{i}\supseteq U_{j}$. The corresponding equivalence classes are also called \df{topological ends} of $X$.  
\end{definition}

To see that the two definitions are indeed equivalent, to each class $[U_{1}\supseteq U_{2}\supseteq\ldots]$ we associate a function $E=E_{[U_{1}\supseteq U_{2}\supseteq\ldots]}$ as in definition \ref{D:Milnor} as follows. Define $E(\partial U_{i}):= U_{i}$ for all $i\in\N$. The function $E$ extends to $\mathcal{K}(X)$. Indeed, given any $K\in\mathcal{K}(X)$ there exists an $i$ such that $K\cap\overline{U_{j}}=\emptyset$ for all $j\geq i$, since $\bigcap\limits_{i\in\N}\overline{U_{i}}=\emptyset$ and $K$ is compact. We define $E(K)$ to be the connected component of $X\setminus K$ containing $U_{i}$. The function $[U_{1}\supseteq U_{2}\supseteq\ldots]\mapsto E_{[U_{1}\supseteq U_{2}\supseteq\ldots]}$ defines a bijection between the set of equivalence classes $[U_{1}\supseteq U_{2}\supseteq\ldots]$ and ${\rm Ends}(X)$ \cite{D}.

\begin{remark}
    \label{R:ends}
From the definitions, one can prove that $\rm X$ has only one end if and only if for every compact set $\rm K\subset X$, there exists a compact set $\rm K'\supset K$ such that $\rm X\setminus K'$ is connected (see \cite[Definition~1.2]{Ray}).\\
\end{remark}
\subsection{Surfaces and the Loch Ness monster}

From now on, \df{surface} will mean a connected \emph{not necessarily compact} \linebreak 2-dimensional real \emph{orientable} manifold with empty boundary unless explicitly specified. A \df{subsurface} of a given  surface is a closed region inside the surface whose boundary consists of a finite number of non-intersecting simple closed curves. The genus of a \emph{compact} bordered surface $S$  with $ q(S)$ boundary curves and Euler characteristic $\chi(S)$ is the number $ g(S)=1-\frac{1}{2}(\chi(S)+q(S))$. A surface is said to be \df{planar} if all of its compact subsurfaces are of genus zero. An end $[U_1\supseteq U_2\supset\ldots]$ is called planar if there exists 
an $i\in\N$ such that $U_i$ is planar.

\begin{definition}
The \df{genus} of a surface is the maximum of the genera of its compact subsurfaces (with boundaries).
\end{definition}
Remark that if a surface $\rm S$ has \df{infinite genus} there exists no finite set $\rm \mathcal{C}$ of mutually non-intersecting simple closed curves with the property that $\rm S\setminus \mathcal{C}$ is  \emph{connected and planar}. We define ${\rm Ends'}(S)\subset{\rm Ends}(S)$ as the set of all ends which are not planar. Two compact surfaces are homeomorphic if and only if they have the same genus. Ker\'ekj\'art\'o's theorem states that two non-compact surfaces $S$ and $S'$ of the same genus and orientability class are homeomorphic if and only if ${\rm Ends'}(S)\subset{\rm Ends}(S)$ and ${\rm Ends'}(S')\subset{\rm Ends}(S')$ are homeomorphic as nested topological spaces. We refer the reader to \cite{R} for an alternative definition of the space of ends (also called \df{the ideal boundary}) in the context of surfaces (compare with \cite{F}) and a proof of Ker\'ekj\'art\'o's theorem.
\begin{definition}
Up to homeomorphism, the \df{Loch Ness monster} is the unique infinite genus surface with only one end.
\end{definition}


\subsection{Ends of groups} In the following paragraphs, we recall facts about ends of groups.

Let $G$ be a finitely generated group and let $A$ be a finite generating set that is closed under inverses. The \textbf{Cayley graph} $\Cayley(G,A)$ of $G$ with respect to $A$ is the graph whose vertex set is $G$ and an edge labeled by $a\in A$ from every vertex $g$ to the vertex $ga$. We regard the Cayley graph of $G$ as a topological space by declaring each edge homeomorphic to the unit interval. We can also turn $\Cayley(G,A)$ into a metric space by declaring each edge \emph{isometric} to the unit interval and considering the path metric. As $A$ is a generating set of $G$, the Cayley graph is connected. The group $G$ acts on the Cayley graph by left multiplication: the edge $g\stackrel{a}{\longrightarrow}ga$ is mapped by the element $h\in G$ to the edge $hg\stackrel{a}{\longrightarrow}hga$. This is a free and discrete action by homeomorphisms (or even isometries) that  respects labels of edges. As $A$ is a finite set, the action of $G$ on its Cayley graph is cocompact.
Indeed, the quotient space $\Cayley(G,A)/G$ is the graph with one vertex and a loop labeled by $a$ for each $a\in A$.


As explicitly stated in the next lemma, the Cayley graph of a group (with respect to any finite basis) captures the topology at infinity of every topological space on which it acts nicely.

\begin{lemma}\label{lem:ends-of-X}
Let $G$ be a finitely generated group, let $A$ be a finitely generating set of $G$ and let $\Cayley(G,A)$ be the Cayley graph of $G$ with respect to $A$.
Let $X$ be a Hausdorff locally compact and locally connected space.
Assume that $G$ acts on $X$ by homeomorphisms properly discontinuously and cocompactly.

Then the space of ends of $\Cayley(G,A)$ and $\Ends(X)$ are homeomorphic.
\end{lemma}

This lemma can be viewed as a topological analog of the celebrated \u{S}varc-Milnor Lemma. A proof can be found in \cite[Corollary 13.5.12]{Geoghegan}.






\begin{corollary}
For any finite generating sets $A$ and $B$ of a group $G$, the spaces of ends $\Ends(\Cayley(G,A))$ and $\Ends(\Cayley(G,B))$ are equal. This is the \textbf{space of ends} $\Ends(G)$ of the group $G$.
\end{corollary}
\proof As $G$ acts freely, discretly and cocompactly on its Cayley graph $X=\Ends(\Cayley(G,B))$, we can apply Lemma~\ref{lem:ends-of-X}. \qed

\medskip

\noindent\textbf{Examples:}

\begin{enumerate}
    \item $\Ends(\Z)=\{-\infty,+\infty\}$.
    \item $\Ends(\Z^2)=\Ends(\Z^3)=\Ends(\Z^n)=\{\infty\}$ for any $n>1$.
    \item If $G$ is a finite group $\Ends(G)=\emptyset$.
    \item The space of ends of a non-trivial finitely generated free group $F_N$, $N>1$ is the Cantor space: it can be identified with the set of infinite reduced words in $A^{\pm 1}$, where $A$ is a basis of $F_N$.
    \item The Coxeter group $[6,4]$ acts properly discontinuously and cocompactly on the hyperbolic plane $\H^2$, thus $\Ends([6,4])=\Ends(\H^2)=\{\infty\}$.
\end{enumerate}

\begin{corollary}\label{cor:ends-subgroup}
Let $G$ be a finitely generated group and $H$ a finite index subgroup of $G$. Then $H$ is finitely generated and $\Ends(G)=\Ends(H)$.
\end{corollary}
\proof The action of $H$ on $X=\Cayley(G,A)$ is free, discrete and cocompact, thus we can apply Lemma~\ref{lem:ends-of-X}.

\medskip

\noindent\textbf{Examples:} \\
\indent 6. The matrix group $SL_2(\Z)$ has a free subgroup of finite index, hence its space of ends is the Cantor space.

\medskip

We emphasize (although we will not use this result) that Stallings completly classified the possible spaces of ends for finitely generated groups:

\begin{theorem}[{\cite{Stallings-1969}}]
The space of ends of a finitely generated group $G$ is one of the followings:
\begin{enumerate}
\item The empty set if and only if $G$ is finite.
\item A set with two elements $\{-\infty;+\infty\}$ if and only if $G$ has a finite index subgroup isomorphic to $\Z$.
\item A Cantor set if and only if $G$ splits as a non-trivial free product or an HNN extension over a finite group.
\item A singleton.
\end{enumerate}
\end{theorem}


\section{Main results}
    \label{S:main}
Consider the set of maps $M(\Gamma,i,\E^2)$ on the plane whose group of symmetries contains two translations with respect to lineraly independent vectors. This set is infinite.
In this section we will show that the surface $S$ on which the minimal regular cover of these maps happens \emph{is always the Loch Ness monster}. This will be achieved in two steps.
\begin{description}
\item[Step 1:] We prove that the surface $S$ has only one end. This follows from lemma \ref{L:1end}, which gives sufficient conditions for a map on a surface to have only one end.
\item [Step 2:] Assuming step 1, we prove that the surface $S$ has infinite genus. For this, suppose first that $S$ has genus. This means that we can find a compact subsurface
\begin{equation}
S'\stackrel{I}{\longhookrightarrow} S
\end{equation}
such that $0<\text{genus}(S')<\infty$. Using the action $\aut (\mathcal{\widetilde{M}})\curvearrowright S$ by homeomorphisms we can translate $S'$ to create an infinite family of embeddings
\begin{equation}
  \label{E:embed}
I_{k}:S'\stackrel{I_{k}}\longhookrightarrow S,\hspace{2mm} k\in\N
\end{equation}
such that $\text{Image}(I_{k})\cap \text{Image}(I_{l})=\emptyset$
if $k\neq l$. This provides infinite genus. Finally, we have to show that $S$ has genus. For this we prove lemma \ref{genus}, which gives sufficient conditions for a map in the plane to have a minimal regular cover with genus. It is easy to see that these conditions are met by any Archimidean map.
\end{description}
\subsection{Step 1} Let $\F(M)$ denote the set of flags of the map $M$ and $\O_{\aut}(M)=\F(M)/\aut(M)$ the corresponding space of $\aut(M)$-orbits. We denote the $\aut(M) $-orbit of a flag $\Phi\in\F(M)$ by $[\Phi]$. Since the action $\mathcal{F}(M)\curvearrowleft\mon(M)$ is transitive, there is a natural (group) morphism
\begin{equation}
  \label{E:kermon}
  \mon(M)\stackrel{f}{\longrightarrow}{\rm Sym}(\O_{\aut}(M)).
\end{equation}
Fix an orbit $[\Phi]$ in $\O_{\aut}(M)$ and define
\begin{equation}
\Stab_{\mon}([\Phi]):=\{r\in\mon(M)\hspace{1mm}|\hspace{1mm} [\Phi r]=[\Phi]\}.
\end{equation}
Clearly,
\begin{equation}
  \label{E:subgroups}
{\rm Ker}(f)\leq\Stab_{\mon}([\Phi])\leq \mon(M).
\end{equation}
\begin{lemma}
    \label{L:kerf}
Let $M(\Gamma,i,\E^2)$ be a map whose group of symmetries contains two translations with respect to lineraly independent vectors.
Then ${\rm Ker}(f)$, as in $(\ref{E:kermon})$, has a finite index subgroup isomorphic to $\mathbb{Z}^n$ for $n\geq 2$.
\end{lemma}
\proof From Proposition \ref{prop:finiteorbits}, we know that $\mathcal{O}_{\aut(M)}$ is finite.
Let $k\in\N$ be the cardinality of this set. Every $k$-tuple $(\Phi_{1},\ldots,\Phi_{k})\in(\mathcal{F}(M))^{k}$ such that $[\Phi_{i}]\neq[\Phi_{j}]$ if $i\neq j$ defines an injective morphism
\begin{equation}
{\rm Ker}(f)\stackrel{g}{\longrightarrow} (\aut(M))^{k}
\end{equation}
in the following way. Recall that every element on $\aut(M)$ is completely determined by the image of one flag in $\mathcal{F}(M)$. Hence, to each $r\in{\rm Ker}(f)$ and $\Phi_{i}$,  $i\in\{1,\ldots, k\}$, we can associate the unique $\alpha_{i}\in\aut(M)$ such that $\Phi_{i}r=\alpha_{i}\Phi_{i}$. The morphism $g$ in injective since the action
$\aut(M)\curvearrowright\F(M)$ is free.
Let $T\cong\Z^{2}$ be the finite index subgroup of $\aut(M)$ consisting of translations and define $A:= g({\rm Ker}(f))\cap(T)^{k}$. Then $g^{-1}(A)$ is a free abelian subgroup which has finite index in ${\rm Ker}(f)$. In what follows we prove that the rank of $g^{-1}(A)$ is bigger than 2. Define the morphism
\begin{equation}
  \label{E:surj}
\Stab_{\mon}([\Phi])\stackrel{h}{\longrightarrow} \aut(M)
\end{equation}
as follows. As seen in the preceding paragraph, for every $r\in\Stab_{\mon}([\Phi])$ there exists a unique $\alpha\in \aut(M)$ such that $\Phi r=\alpha \Phi$. We claim that $h$ is surjective. Indeed, recall that $\mon(M)$ always acts transitively on $\F(M)$, hence for every $\alpha\in\aut(M)$ there exists $r\in\mon(M)$ such that $\Phi r=\alpha \Phi$. Since $h$ is surjective, $h^{-1}(T)$ is a free abelian subgroup of finite index in $\Stab_{\mon}([\Phi])$ of rank at least 2. Hence:
\begin{equation}
G:=h^{-1}(T)\cap g^{-1}(A)\leq{\rm Ker}(f)
\end{equation}
is a free abelian subgroup of finite index and rank at least 2. \qed

\begin{lemma}
    \label{L:1end}
Let $M=M(\Gamma,i,S)$ be a map such that $\aut(M)$ contains a finite index subgroup isomorphic to $\Z^{n}$ with $n\geq 2$ and $\mathcal{O}_{\aut}(M)$ is a finite set. Let $\widetilde{M}=\widetilde{M}(\widetilde{\Gamma},\tilde{i},\widetilde{S})$ be the minimal regular cover of $M$. Then ${\rm Ends}(\widetilde{S})$ is just one point.
\end{lemma}
\proof For any map $M=M(\Gamma,i,S)$, the minimal regular cover $\widetilde{M}=\widetilde{M}(\widetilde{\Gamma},\tilde{i},\widetilde{S})$ satisfies that:
\begin{equation}
  \label{E:AutMon}
\mon(\widetilde{M})\simeq\aut(\widetilde{M})\simeq\mon(M).
\end{equation}
Since $\widetilde{M}$ is regular, the action by homeomorphisms  $\aut (\widetilde M)\curvearrowright\widetilde S$ is properly discontinuous and cocompact (see $\S \ref{S:tess}$). Hence, by lemma \ref{lem:ends-of-X}, the spaces ${\rm Ends}(\widetilde S)$ and ${\rm Ends} (\aut(\widetilde{M}))$ are homeomorphic. Since $\mathcal{O}_{\aut}(M)$ is finite, ${\rm Ker}(f)$ has finite index in $\mon (M)$. From lemma \ref{L:kerf} we know that ${\rm Ker}(f)$ contains a finite index subgroup isomorphic to $\mathbb{Z}^n$ for some $n\geq 2$. This concludes the proof of the lemma. \qed

\subsection{Step 2}
\begin{lemma}
  \label{genus}
Let $M(\Gamma,i,\E^{2})$ be a map satisfying the hypotheses of lemma \ref{egamma}, then the surface  $\widetilde{S}$ corresponding to the minimal regular cover  has genus.
\end{lemma}
\proof The proof is by contradiction. Let $\widetilde{S}\stackrel{\pi}{\longrightarrow} \E^{2}$ be the (ramified) covering map and suppose that $\widetilde{S}$ has no genus. Since it has only one end, then $\widetilde{S}$ is homeomorphic to $\E^{2}$. Now consider $\gamma$ as in lemma \ref{egamma} and let $\pi^{-1}(\gamma)=\{\widetilde{\gamma_{j}}\}_{j\in\N}$, where $\gamma_{j}$ is a simple closed curve for all $j$. Remark that by the same lemma, the curve $\gamma$ does not contain branch points of $\pi$. Since $\gamma$ is a Jordan curve, it defines a bounded component $B$ of $\E^{2}$. Every simple closed curve $\widetilde{\gamma_{j}}$ defines a bounded component $\widetilde{B_{j}}$ of $\widetilde{S}$. Given that the number of flags of the minimal regular cover contained in $\widetilde{B_{j}}$ must be constant, we can assume that $\widetilde{\gamma_{j}}\cap\widetilde{B_{k}}=\emptyset$ if $j\neq k$. That is, the curves forming $\pi^{-1}(\gamma)$ are not nested.\\
\indent The map $M$ has at least two faces of different size and $\aut(M)$ contains two linearly independent translations, hence there are infinitely many branch points of $\pi$ in $\E^{2}$ arranged like points in a lattice. Hence we can pick two branch points $x\neq y$ of $\pi$ in $\E^{2}$, such that $x\in B$ and $y\in\E^{2}\setminus B$. We consider two cases, both leading to a contradiction.
\begin{description}
\item [Case A]. There exist a ramification point $p\in\pi^{-1}(x)$ such that $p\in \widetilde{B_{j}}$. This situation leads to a contradiction as follows.  Consider  an arc $I\subset\E^{2}$ joining  the branch point $x$ to a point $x'\in\gamma$ and containing no branch points of $\pi$. We can  pick two lifts $\widetilde{I}$, $\tilde{I}'$ of $I$ with an endpoint at $p$.
Since $I$ has no branch points:
$$
\{\widetilde{I}\setminus\{p,\pi^{-1}(x')\}\}\cap\{\widetilde{I}'\setminus\{p,\pi^{-1}(x')\}\}=\emptyset.
$$
On the other hand, $\pi$ restricted to any $\widetilde{\gamma_{j}}$ has to be injective, hence both $\tilde{I}$ and $\widetilde{I}'$ join $p$ with a point $p'\in\widetilde{\gamma_{j}}$. But this implies that $\pi(p')$ is a branch point in $\gamma$, which is a contradiction since we supposed that $\gamma$ contains no branch points.
\item [Case B]. For all $j\in\N$, $\pi^{-1}(x)$ and $\pi^{-1}(y)$ are disjoint from $\widetilde{B_{j}}$. Pick $p\in\pi^{-1}(x)$ and $q\in\pi^{-1}(y)$. Since $\pi^{-1}(\gamma)$ is closed and $\pi$ is a branched covering, we can find a path $\eta$ joining $p$ to $q$ and contained in $\widetilde{S}\setminus \bigcup\limits_{j\in\N} \widetilde{B_{j}}$. But then $\pi(\eta)$ would be a path joining $x$ to $y$ and not intersecting $\gamma$, which is a contradiction.\qed
\end{description}


From the preceeding lemmas we can conclude the main result of this article.
\begin{theorem}
    \label{archlochness}
Let $M=M(\Gamma,i,\E^2)$ be a map with convex faces with at least two faces different sizes such that it symmetry group contains two translations with respecto to linearly independent vectors. Then the surface $S$ on which the minimal regular cover of $M$ takes place is homeomorphic to the Loch Ness Monster.
\end{theorem}
\begin{remark}
There is a dual statement for Theorem \ref{archlochness} where we ask the map $M=M(\Gamma,i,\E^2)$ to have convex faces with 2 vertices of different degrees and we make no requirements on the size of the faces of $M$.
\end{remark}
Since all Archimedean maps satisfy the hypothesis of Theorem  \ref{archlochness} we have the following:
\begin{corollary}
    \label{C:arch}
The minimal regular cover of any Archimedean map takes place on a surface homeomorphic to the Loch Ness monster.
\end{corollary}

\section{An example}

In this section we present a more concrete proof of a particular case of Theorem \ref{archlochness}, namely that the surface on which the minimal regular cover of the Archimedean tiling 3.6.3.6 lies is the Loch Ness monster. 

We prove first that the surface has only one end, which is done by considering the dual graph of the tiling, and realizing it is the Cayley graph of a group with only one end. Then we proceed to prove the surface has infinite genus using Euler's characteristic. 

We start with a description on the minimum regular cover of 3.6.3.6, as first described in \cite{GorDan}. Recall that $\{6,4\}$ is the hyperbolic tiling of $\H$ having four hexagons around each vertex:
\begin{center}
    \includegraphics[scale=0.5]{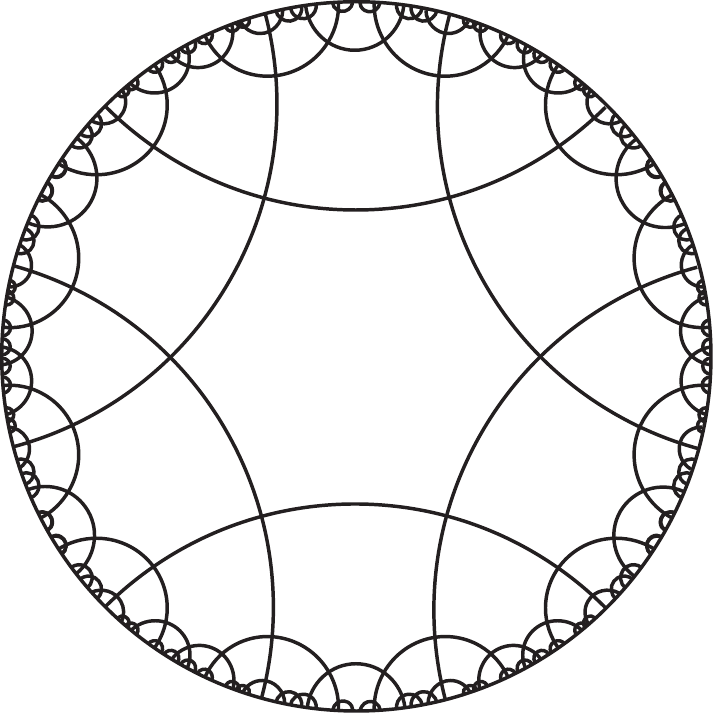}\\
    \small{Figure 4: The $\{6,4\}$ hyperbolic tiling.}
\end{center}

Then, the minimal regular cover can be constructed by applying to the tiling $\{6,4\}$, whenever possible and in every direction, the two identifications $A$ and $B$ in Figures 5 and 6 respectively.
\begin{center}
    \includegraphics[scale=0.5]{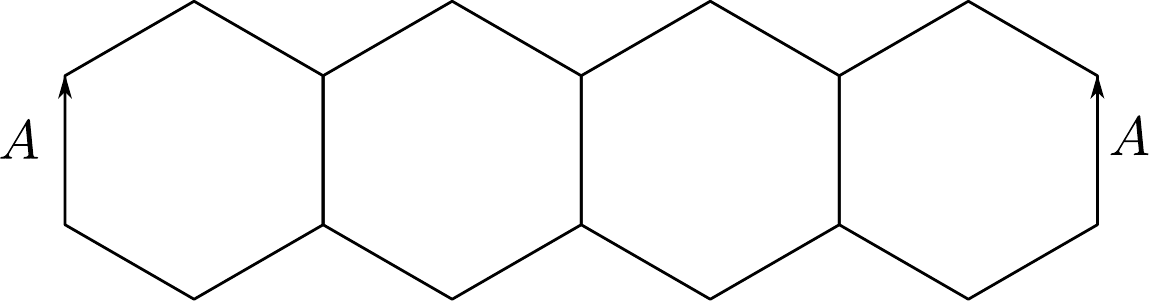} \\
    \small{Figure 5: Identification A, $[(\rho_1 \rho_0)^2 \rho_1 \rho_2]^4$.}
\end{center}

\begin{center}
    \includegraphics[scale=0.5]{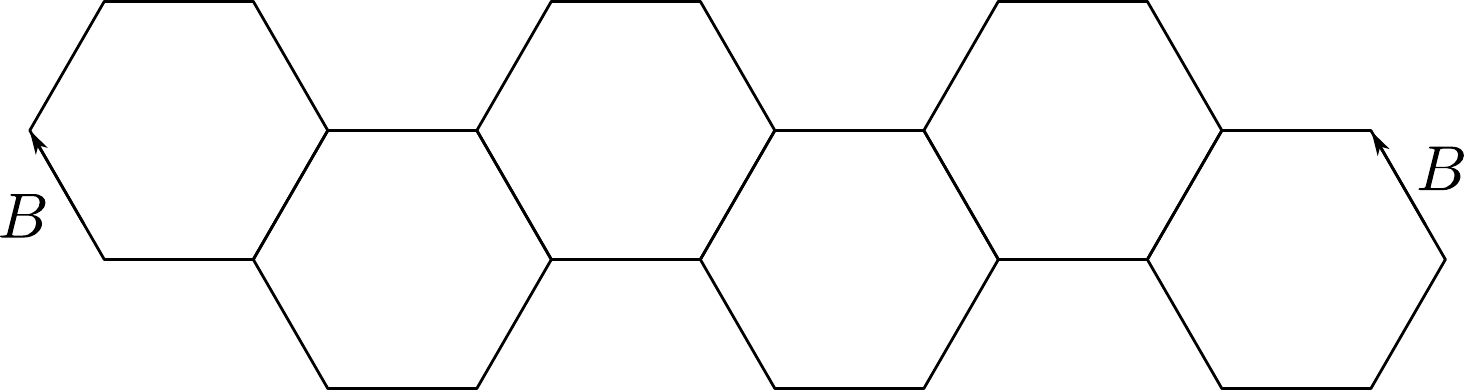} \\
    \small{Figure 6: Identification B, $[(\rho_1 \rho_0)^2 \rho_2]^6$.}
\end{center}

The minimum regular cover is then realized on the surface
\begin{equation}
S:=\H/\sim,
\end{equation}
where $\sim$ is the equivalence relation defined by $A$ and $B$.

Let $\cD$ be the set of hexagons that define the tiling of $S$. Construct the dual graph $\Lambda$ of the map. That is, construct the graph whose vertices are $\cD$ and there is an edge between two hexagons if and only if they share an edge. Note that any finite set of $V(\Lambda)$ corresponds to a compact set of $S$ (by taking the closed hexagons). Also, any connected subgraph of $\Lambda$ corresponds to a connected subsurface of $S$.

\begin{lemma}\label{Soneend}
Let $\Lambda$ be the dual graph of the minimum regular cover of 3.6.3.6. Then the edges of $\Lambda$ can be labeled in such a way that $\Lambda$ is the Cayley graph of the group $G \cong H \times H$ with
    $$H:=\langle a,b,c\ :\ a^{2}=b^{2}=c^{2}=(abc)^{2}=1\rangle.$$
\end{lemma}

\proof 
We note that every vertex of $\Lambda$ has degree 6. We shall color the edges with 6 colors ($\{1,2,3,4,5,6\}$) in such a way that if two edges share a vertex then the colors are different. In particular, we prove that the chromatic number of the line graph is exactly 6. The edges of the dual graph are in an obvious one-to-one correspondence to the edges of the original map.

Pick any hexagon of $S$ and call it the base hexagon. Color its incoming edges in a clockwise order with colors 1 through 6. Then color the adjacent hexagons in a counter-clockwise fashion.
  \begin{center}
    \includegraphics[scale=0.4]{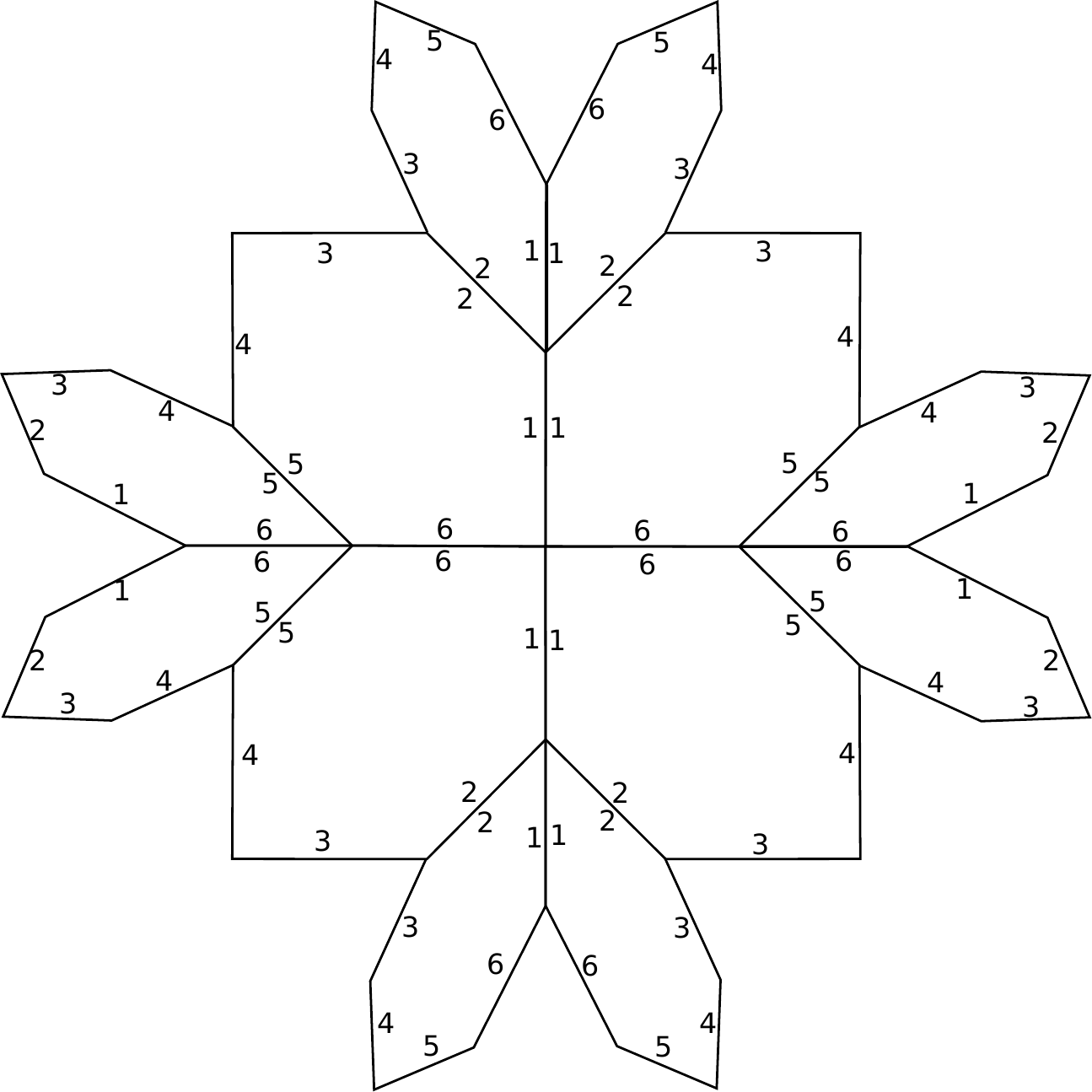} \\
    \small{Figure 7: Twelve hexagons with their coloring.}
  \end{center}

Observe that at every vertex of the original map, opposite edges have the same color. Observe also that identifications $A$ and $B$ are consistent with the coloring (\emph{i.e.}, the color of the two identified edges match and so does the orientation in both cases).

We think of the symbol $a_i$ as the action of walking by the edge with color $i$ in $\Lambda$. Consider the free group $F_6$ with generators $\{a_i:i\in \Z_6\}$. Then we can associate to each word $a_{i_1}a_{i_2}...a_{i_n}\in F_6$ a hexagon, namely the hexagon that results from walking from the hexagon marked as the \emph{base hexagon} first by the edge of color $i_1$, then color $i_2$ and so on. 

Of course, there may be many ways to reach each hexagon and so we may consider a quotient of $F_6$. We have the following relations:
\begin{enumerate}
    \item $a_{i}^2 = e$, since walking the same edge twice takes you back to the starting hexagon.
    \item $a_{i}a_{i+1} = a_{i+1}a_{i}$, since there are two ways of getting to a hexagon that shares only a vertex with the hexagon you are currently at.
    \item $a_{i}a_{i+3} = a_{i+3}a_{i}$, because of identification $A$.
    \item $(a_{1}a_{3}a_{5})^{2}=(a_{2}a_{4}a_{6})^{2}=e$, because of identification $B$.
\end{enumerate}

These are all the required relations. The first two come from the hyperbolic tiling $\{6,4\}$ and the last two are implied by $A$ and $B$.

We conclude that $\Lambda$ is the Cayley graph of the group $G$ on 6 generators $a_{1},\ldots,a_6$ with only the following relations:
\begin{equation}
a_{i}^{2}=e,\ a_{2i}a_{2j+1}=a_{2j+1}a_{2i},\ (a_{1}a_{3}a_{5})^{2}=(a_{2}a_{4}a_{6})^{2}=e.
\end{equation}

Note that since $a_{i}$ and $a_{j}$ commute whenever $i$ and $j$ have opposite parity, we can write $G \cong H\times H$ with
$$H:=\langle a,b,c\ :\ a^{2}=b^{2}=c^{2}=(abc)^{2}=1\rangle.$$ \qed

We now prove that $H$ has only one end.

\begin{lemma}\label{deletefiniteH}
    Let $H=\langle a,b,c\ :\ a^{2}=b^{2}=c^{2}=(abc)^{2}=1\rangle$ and let $F \subset V(\Cayley(H,\{a,b,c\}))=H$ be a finite set of vertices of the graph. Then there exists $F'$ finite with $F\subset F' \subset H$ such that the graph induced on $H\setminus F'$ is connected.
\end{lemma}

\proof We claim that the Cayley graph of $H$ with respect to $a$, $b$ and $c$ is the 1-skeleton of the tiling by hexagons of the plane:

\begin{center}
    \includegraphics[scale=0.5]{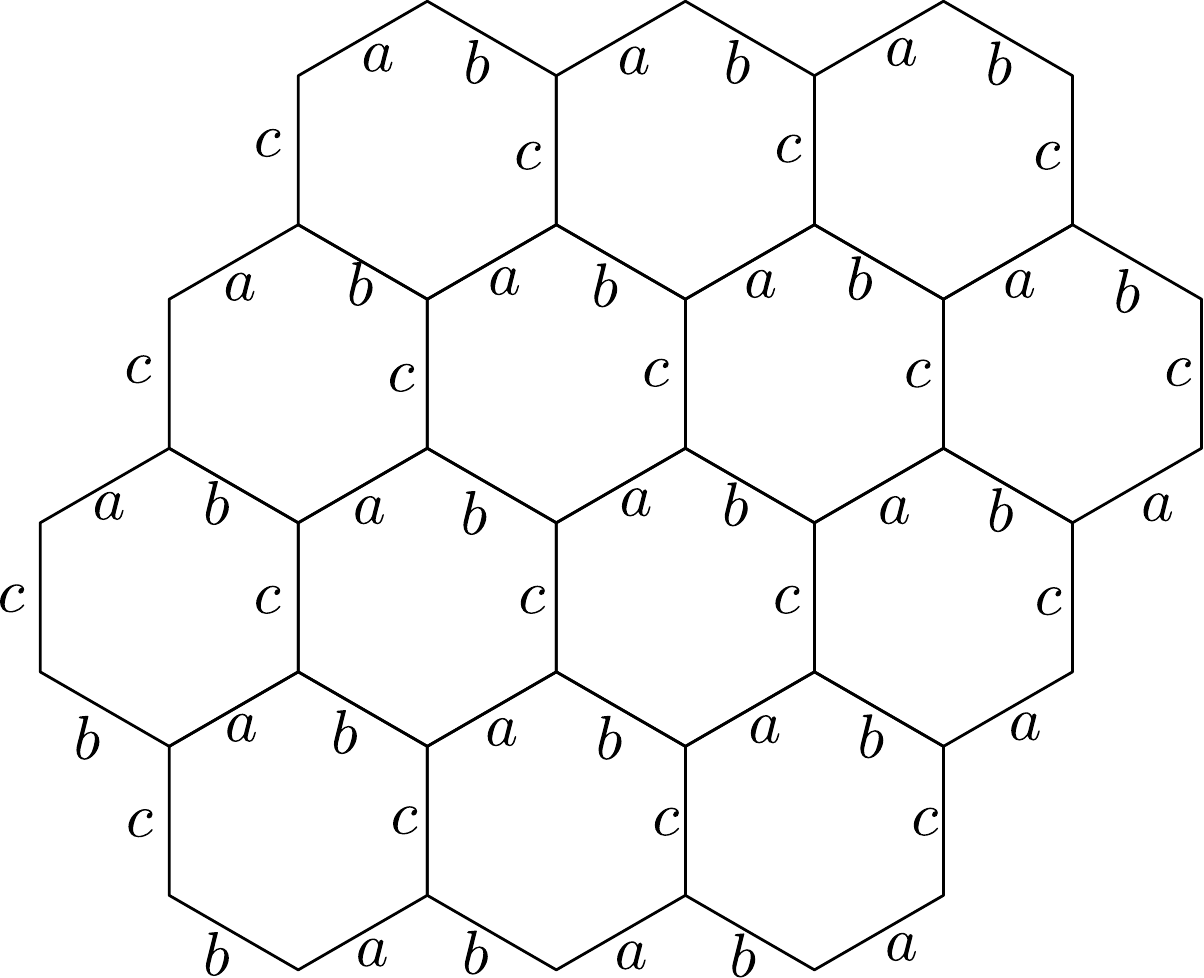}\\
    \small{Figure 8: The Cayley graph of $H$.}
\end{center}

First of all, note that  $deg(h)=3$ for all $h \in V(\Cayley(H, \{a,b,c\}))$.
Then, the following cycles
\[
\gamma_{1}: 1, a, ab, abc, cb, c, 1,
\]
\[
\gamma_{2}: 1, b, bc, bca, ac, a, 1,
\]
\[
\gamma_{3}: 1, c, ca, cab, ba, b, 1
\]
satisty that $\gamma_{i} \subset Cayley(H,\{a,b,c\})$ for $i \in \{1,2,3 \}$, $\gamma_{1}\cap \gamma_{2}\cap \gamma_{3}=1$ (the vertex corresponding to 1), $\gamma_{i}\cap\gamma_{j}$ is only one
edge $\forall \ i\neq j \in \{1,2,3 \}$ and each vertex $v \in \gamma_{i}$ is adjacent with only two vertices in $\gamma_{i}$.
We know  that every Cayley graph is vertex-transitive, \emph{i.e.}, the above happens on each vertex of $\Cayley(H,\{a,b,c\})$.
Notice that each edge in $E(\Cayley(H,\{a,b,c\}))$ is in only two hexagons of this type, since the degree at every vertex is 3. Also notice that there are no other relations than the ones implied by the tiling of $\E^2$ by hexagons and so $\Cayley(H,\{a,b,c\})$ must be this tiling.

Then the lemma becomes obvious: If $F$ is a finite set of the vertices, we can take a set of vertices $F'$ covered by a ball that also covers all of $F$. Then the complement of $F'$ is connected. \qed

Note that this lemma implies that the same will be true for
    $$\Cayley(G,\{a_1,a_2,a_3,a_4,a_5,a_6\}),$$
since $G \cong H\times H$.
\begin{lemma}\label{OneEnd} The surface $S$, on which the minimal regular cover of the \linebreak Archimedean tiling 3.6.3.6 lies, has only one end. \end{lemma}
\proof Let $K$ be a compact set of $S$. Suppose it intersects a set $F$ of hexagons. Applying lemmas \ref{Soneend} and \ref{deletefiniteH} we conclude that there only remains a single infinite connected component of $S \setminus K$. Then, by lemma~\ref{lem:ends-of-X} we conclude that $S$ has only one end.\qed
\begin{lemma}\label{Sinfinitegenus} The surface $S$, on which the minimal regular cover of the \linebreak Archimedean tiling 3.6.3.6 lies, has infinite genus. \end{lemma}

\proof We proceed by contradiction. From the previous lemma we know that $S$ has only one end.

Either the genus $g$ of the surface $S$ is in $\{0,1,2,...\}$, or $S$ has infinite genus. Notice that if $g \geq 1$, then, by the argument given on Step 2 at the beginning of Section \ref{S:main}, automatically $S$ has infinite genus. So, let us assume for the sake of contradiction that $g=0$, which would mean that surface $S$ is homeomorphic to $\E^2$.

Consider 4 hexagons as in identification $A$:

\begin{center}
    \includegraphics[scale=0.4]{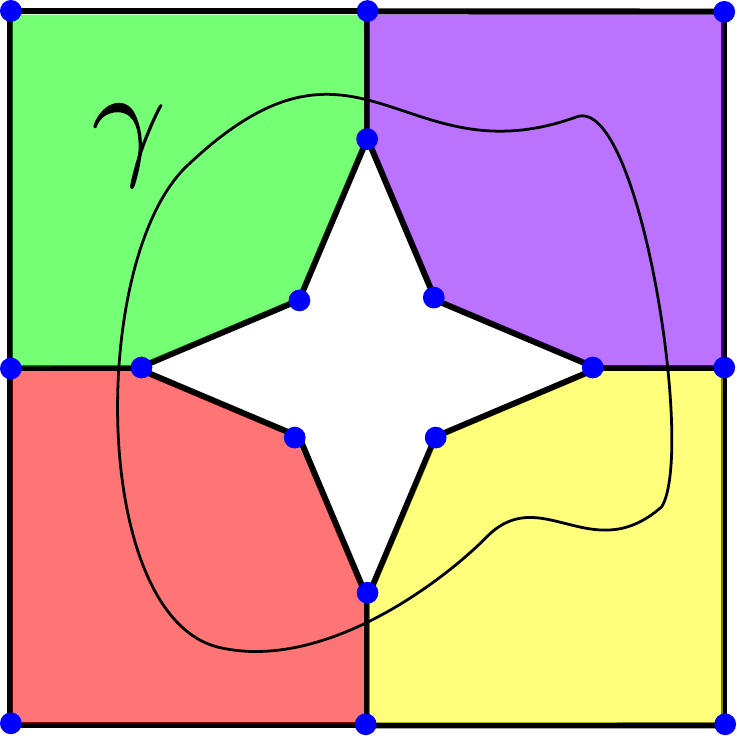} \\
    \small{Figure 9: Four hexagons in the plane, with identification $A$.}
\end{center}

Then we consider the graph $\hat{H}$ built from the four hexagons in the figure and all the hexagons inside the inner region (dismiss all hexagons in the outer region). 
Note that if there were infinitely many hexagons in the inner region then the set of vertices would have an accumulation point, contradicting the definition of map given in section 2. So, we may assume there are only a finite number of hexagons in the inner region.

On the other hand, $\hat{H}$ is embedded in $\E^2$. We remark on this graph:

\begin{itemize}
    \item There are four vertices of degree 3, four vertices of degree 2, and all the others (say $v$) of degree 4.
    \item There is one octagon (the outer octagon), \emph{i.e.}, one face with only eight edges, and all the others (say $h$ of them) are hexagons.
    \item For the edges, we count them in two ways.
        $$E(\hat H)=\frac{6h+8}{2}=\frac{4v+8+12}{2}.$$
\end{itemize}

Now, appliying Euler's characteristic, we get:

$$2=V(\hat H)-E(\hat H) + F(\hat H) = (v+8) - (2v+4+6) + (\frac{2v}{3} + 2 + 1) = \frac{-v}{3}+1.$$

Since $v$ is positive, this yields a contradiction. \qed

\begin{theorem}
The surface S on which the minimal regular cover of the \linebreak Archimedean tiling 3.6.3.6 lies; is homeomorphic to the Loch Ness Monster.
\end{theorem}

\proof 
 Lemma \ref{OneEnd} together with lemma \ref{Sinfinitegenus} yield that $S$ is the \emph{Loch Ness Monster}. \qed

\vspace{0.3cm}

We believe the same kind of concrete and arguably more illustrative proof will work for the other seven Archimedean tilings.


\section{Final remarks}

The minimal regular cover of a map satisfying the hypotheses of Theorem \ref{archlochness} occurs on a surface $\widetilde{S}$ that has a natural hyperbolic metric. In this section we describe how this metric arises and state some questions concerning the corresponding Teichm\"uller spaces.\\
\indent Let $M=M(\Gamma,i,\E^2)$ be a map with convex faces with at least two faces of different sizes such that its symmetry group contains two translations with respect to linearly independent vectors. Let
\begin{equation}
\begin{array}{ccc}
p & = & LCM\{p_i\hspace{1mm}|\hspace{1mm} \text{$p_i=$\# of edges in a face of $M$}\}\\
\\
q & = & LCM\{q_i\hspace{1mm}|\hspace{1mm} \text{$q_i$=degree of a vertex in $M$}\}
\end{array}
\end{equation}
The map $M$ is covered by the (regular) tessellation $\{p,q\}$ of $\H$. Recall from \S \ref{S:tess} that the automorphism and symmetry group of $\{p,q\}$ are isomorphic to the Coxeter group $[p,q]$. Since the tessellation is regular,  $[p,q]$ is also isomorphic to the monodromy group of $\{p,q\}$. Let $\Phi$ be a flag in  the minimal cover $\widetilde{M}$ and denote by $Stab_{[p,q]}(\Phi)$ the stabilizer of the flag $\Phi$ in the monodromy group.
Then, the minimal regular covering of $M$ is given by
\begin{equation}
    \label{E:teich}
\widetilde{M}=\{p,q\}/Stab_{[p,q]}(\Phi).
\end{equation}
Since $\widetilde{M}$ is the minimal cover of $M$ we have that $\mon(\widetilde{M})\cong\mon(M)$. From \cite{GorDan} we know that
 $Stab_{[p,q]}(\Phi)=\bigcap\limits_{{\Psi}\in\mathcal{F}(M)}Stab_{[p,q]}(\psi)$. In other words, we can express the minimal regular cover of $M$ as a quotient of the hyperbolic tilling $\{p,q\}$ by a subgroup of isometries of $\H$.\\
\indent Unlike the case of surfaces of finite topological type, there are several different Teichm\"uller spaces that are associated to a surface of topologically infinite type (see \cite{Liu} for a precise statement). From Theorem \ref{archlochness} we know that the hyperbolic surfaces (\ref{E:teich}) have all the same topological type.
\begin{question}
Do all the hyperboic surfaces defined by \ref{E:teich} lie in the same Teichm\"uller space?
\end{question}


\begin{bibdiv}
  \begin{biblist}
 
\bib{conderdob}{article}{
AUTHOR = {Conder, Marston},
AUTHOR = {Dobcs\'anyi, Peter},
TITLE = {Determination of all regular maps of small genus},
JOURNAL = {J. Combin. Theory Ser. B},
FJOURNAL = {81},
VOLUME = {2},
YEAR = {2001},
PAGES = {224--242}}


\bib{Con}{book}{
    AUTHOR = {Conway, John},
    AUTHOR = {Burgiel, Heidi},
    AUTHOR = {Goodman-Strauss, Chaim},
     TITLE = {The symmetries of things},
      YEAR = {2008},
 PUBLISHER = {A K Peters, Wellesley},
     }

 \bib{D}{article}{
 AUTHOR = {Diestel, Reinhard},
 AUTHOR =  {K{\"u}hn, Daniela},
     TITLE = {Graph-theoretical versus topological ends of graphs},
      NOTE = {Dedicated to Crispin St. J. A. Nash-Williams},
   JOURNAL = {J. Combin. Theory Ser. B},
  FJOURNAL = {Journal of Combinatorial Theory. Series B},
    VOLUME = {87},
      YEAR = {2003},
    NUMBER = {1},
     PAGES = {197--206},
 }
 
 
 
\bib{F}{article}{
AUTHOR = {Freudenthal, Hans},
     TITLE = {\"{U}ber die {E}nden topologischer {R}\"aume und {G}ruppen},
   JOURNAL = {Math. Z.},
  FJOURNAL = {Mathematische Zeitschrift},
    VOLUME = {33},
      YEAR = {1931},
    NUMBER = {1},
     PAGES = {692--713},
         }

\bib{Geoghegan}{book}{
    AUTHOR = {Geoghegan, Ross},
     TITLE = {Topological methods in group theory},
    SERIES = {Graduate Texts in Mathematics},
    VOLUME = {243},
 PUBLISHER = {Springer},
   ADDRESS = {New York},
      YEAR = {2008},
     PAGES = {xiv+473},
      ISBN = {978-0-387-74611-1},
   MRCLASS = {57M07 (20F65 20J05 55-02 55P57 57-02)},
  MRNUMBER = {2365352 (2008j:57002)},
MRREVIEWER = {John G. Ratcliffe},
       DOI = {10.1007/978-0-387-74614-2},
       URL = {http://dx.doi.org/10.1007/978-0-387-74614-2},
}


\bib{Ghys}{article}{
AUTHOR = {Ghys, {\'E}tienne},
     TITLE = {Topologie des feuilles g\'en\'eriques},
   JOURNAL = {Ann. of Math. (2)},
  FJOURNAL = {Annals of Mathematics. Second Series},
    VOLUME = {141},
      YEAR = {1995},
    NUMBER = {2},
     PAGES = {387--422},
     }

\bib{GruShe}{book}{
AUTHOR = {Gr\"unbaum, Branko}
AUTHOR =  { Sheppard, Geoffrey Colin},
     TITLE = {Tilings and Patterns},
      YEAR = {1987},
 PUBLISHER = {Freeman, New York},
     }


\bib{Hartley}{article}{
AUTHOR = {Hartley, Michael},
     TITLE = {All polytopes are quotients, and isomorphic polytopes are quotients by conjugate subgroups},
   JOURNAL = {Discrete Comput. Geom.},
  FJOURNAL = {21},
    VOLUME = {2},
      YEAR = {1999},
     PAGES = {289--298},
     }

\bib{PelGor}{article}{
AUTHOR = {Hartley, Michael}
AUTHOR = {Pellicer, Daniel}
AUTHOR = {Williams, Gordon},
     TITLE = {Minimal covers of the prisms and antiprisms},
   JOURNAL = {Discrete Math.},
  FJOURNAL = {312},
    VOLUME = {20},
      YEAR = {2012},
     PAGES = {3046--3058},
     }

\bib{HarWil}{article}{
AUTHOR = {Hartley, Michael},
AUTHOR = {Williams, Gordon},
     TITLE = {Representing the sporadic Archimedean polyhedra as abstract polytopes},
   JOURNAL = {Discrete Math.},
  FJOURNAL = {310},
    VOLUME = {12},
      YEAR = {2010},
     PAGES = {1835--1844},
     }

\bib{Gareth}{article}{
AUTHOR = {Jones, Gareth A.},
     TITLE = {Regular embeddings of complete bipartite graphs: classification and enumeration},
   JOURNAL = {Proc. Lond. Math. Soc.},
  FJOURNAL = {101},
    VOLUME = {2},
      YEAR = {2010},
     PAGES = {427--453},
     }

\bib{L}{article}{
AUTHOR = {Lefschetz, S.},
     TITLE = {Book {R}eview: {V}orlesungen \"uber {T}opologie. {I}.
              {F}l\"achentopologie},
   JOURNAL = {Bull. Amer. Math. Soc.},
  FJOURNAL = {Bulletin of the American Mathematical Society},
    VOLUME = {31},
      YEAR = {1925},
    NUMBER = {3-4},
     PAGES = {176},
}

\bib{Liu}{article}{
AUTHOR = {Liu, Lixin},
AUTHOR = {Papadopoulos, Athanase},
     TITLE = {Some metrics on {T}eichm\"uller spaces of surfaces of infinite
              type},
   JOURNAL = {Trans. Amer. Math. Soc.},
  FJOURNAL = {Transactions of the American Mathematical Society},
    VOLUME = {363},
      YEAR = {2011},
    NUMBER = {8},
     PAGES = {4109--4134},
}

\bib{ARP}{book}{
AUTHOR = {McMullen, Peter},
AUTHOR = { Schulte, Egon},
     TITLE = {Abstract Regular Polytopes},
      YEAR = {2002},
PUBLISHER = {Cambridge University Press},
     }


\bib{Mil}{article}{
AUTHOR = {Milnor, John W.},
     TITLE = {Infinite cyclic coverings},
 BOOKTITLE = {Conference on the {T}opology of {M}anifolds ({M}ichigan
              {S}tate {U}niv., {E}. {L}ansing, {M}ich., 1967)},
     PAGES = {115--133},
 PUBLISHER = {Prindle, Weber \& Schmidt, Boston, Mass.},
      YEAR = {1968},
}

\bib{M}{article}{
AUTHOR = {Mazurkiewicz, Stefan},
     TITLE = {Recherches sur la th\'eorie des bouts premiers},
   JOURNAL = {Fund. Math.},
  FJOURNAL = {Polska Akademia Nauk. Fundamenta Mathematicae},
    VOLUME = {33},
      YEAR = {1945},
     PAGES = {177--228},
}


\bib{Mixer}{article}{
AUTHOR = {Mixer, Mark},
AUTHOR = {Pellicer, Daniel},
AUTHOR = {Williams, Gordon},
     TITLE = {Minimal Covers of the Archimedean Tilings, Part II},
   JOURNAL = {In preparation}
      YEAR = {2012},
     PAGES = {10},
     }

\bib{Monson}{article}{
AUTHOR = {Monson, Barry},
AUTHOR = {Pellicer, Daniel},
AUTHOR = {Williams, Gordon},
     TITLE = {Mixing and monodromy of abstract polytopes},
   JOURNAL = {To appear in Trans. Amer. Mat. Soc.},
    PAGES = {10},
     }

\bib{Barry}{article}{
AUTHOR = {Monson, Barry},
AUTHOR = {Pellicer, Daniel},
AUTHOR = {Williams, Gordon},
 TITLE = {The Tomotope},
 JOURNAL = {To appear in Ars Math. Contemp.},
  PAGES = {10},
     }



\bib{GorDan}{article}{
AUTHOR = {Pellicer, Daniel}
AUTHOR = {Williams, Gordon},
     TITLE = {Minimal Covers of the Archimedean Tilings, Part I},
   JOURNAL = {To appear in Electron. J. Combin.},
     PAGES = {10},
     }

\bib{PSV}{article}{
AUTHOR = {Przytycki, Piotr},
AUTHOR = {Schmith{\"u}sen, Gabriela},
AUTHOR = {Valdez,Ferr{\'a}n},
     TITLE = {Veech groups of {L}och {N}ess monsters},
   JOURNAL = {Ann. Inst. Fourier (Grenoble)},
  FJOURNAL = {Universit\'e de Grenoble. Annales de l'Institut Fourier},
    VOLUME = {61},
      YEAR = {2011},
    NUMBER = {2},
     PAGES = {673--687},
}

\bib{Ray}{article}{
AUTHOR = {Raymond, Frank},
     TITLE = {The end point compactification of manifolds},
   JOURNAL = {Pacific J. Math.},
  FJOURNAL = {Pacific Journal of Mathematics},
    VOLUME = {10},
      YEAR = {1960},
     PAGES = {947--963},
     }

\bib{R}{article}{
AUTHOR = {Richards, Ian},
     TITLE = {On the classification of noncompact surfaces},
   JOURNAL = {Trans. Amer. Math. Soc.},
  FJOURNAL = {Transactions of the American Mathematical Society},
    VOLUME = {106},
      YEAR = {1963},
     PAGES = {259--269},
     }

\bib{Ratcliffe}{book}{
    AUTHOR = {Ratcliffe, John G.},
     TITLE = {Foundations of hperbolic manifolds},
      YEAR = {1994},
 PUBLISHER = {Springer-Verlag New York},
     }

\bib{Stallings-1969}{book}{
AUTHOR = {Stallings, John},
TITLE = {Group theory and three-dimensional manifolds},
YEAR={1969},
PUBLISHER= {Yale University Press, New Haven, Conn.-London, 1971},
}

\bib{V1}{article}{
    AUTHOR = {Valdez, Ferr{\'a}n},
     TITLE = {Billiards in polygons and homogeneous foliations on $\mathbf{C^2}$},
   JOURNAL = {Ergodic Theory Dynam. Systems},
  FJOURNAL = {Ergodic Theory and Dynamical Systems},
    VOLUME = {29},
      YEAR = {2009},
    NUMBER = {1},
     PAGES = {255--271},
}

\bib{V2}{article}{
    AUTHOR = {Valdez, Ferr{\'a}n},
     TITLE = {Infinite genus surfaces and irrational polygonal billiards},
   JOURNAL = {Geom. Dedicata},
  FJOURNAL = {Geometriae Dedicata},
    VOLUME = {143},
      YEAR = {2009},
     PAGES = {143--154},
      ISSN = {0046-5755},
}

\bib{Wilson}{article}{
    AUTHOR = {Wilson, Steve},
     TITLE = {Families of regular graphs in regular maps},
   JOURNAL = {J. Combin. Theory Ser. B},
  FJOURNAL = {Journal of Combinatorial Theory Series B},
    VOLUME = {85},
      YEAR = {2002},
     PAGES = {269--289},
      ISSN = {0046-5755},
}



     
\end{biblist}
  \end{bibdiv}
\end{document}